%%To be processed in msamsppt.sty

\input amstex
\documentstyle{msamsppt}
\pageno=1
\magnification1200
\catcode`\@=11
\def\logo@{}
\catcode`\@=\active
\NoBlackBoxes
\vsize=23.5truecm
\hsize=16.5truecm

\let\Bbb=\sanssymbol
\let\qed=\relax
\let\square=\relax

\def\d{d\!@!@!@!@!@!{}^{@!@!\text{\rm--}}\!}

\def\leg{\;\dot{\le}\;}

\def\eg{\;\dot {=}\;}
\def\crp{\overline{\Bbb R}_+}
\def\crm{\overline{\Bbb R}_-}

\def\rnp{{\Bbb R}^n_+}

\def\rnpm{\Bbb R^n_\pm}

\def\ang#1{\langle {#1} \rangle}

\def\rp{ \Bbb R_+}
\def\rmi{ \Bbb R_-}

\define\tr{\operatorname{tr}}
\define\Tr{\operatorname{Tr}}
\define\res{\operatorname{res}}

\define\OP{\operatorname{OP}}
\define\deln{\tfrac{\partial  _\lambda ^{N-1}}{(N-1)!}}
\define\Res{\operatornamewithlimits{Res}}

\document

\topmatter
\title Logarithms and sectorial projections for elliptic
 boundary problems
\endtitle
\author
Anders Gaarde and Gerd Grubb
\endauthor
\affil
{Copenhagen Univ\. Math\. Dept\.,
Universitetsparken 5, DK-2100 Copenhagen, Denmark.
E-mail {\tt gaarde\@math.ku.dk} and {\tt grubb\@math.ku.dk}}\endaffil
\rightheadtext{Logarithms and sectorial projections}
\abstract On a compact manifold with boundary, consider the
realization $B$ of an elliptic, possibly pseudodifferential, boundary
value problem having a spectral cut (a ray free of eigenvalues), say
${\Bbb R}_-$. In the first part of the paper we define and discuss 
in detail the operator $\log B$;
its residue (generalizing the Wodzicki residue) is essentially 
proportional to the zeta function value at zero, $\zeta (B,0)$, and
it enters in an important way in studies of composed zeta functions $\zeta
(A,B,s)=\Tr(AB^{-s})$ (pursued elsewhere).

There is a similar definition of the operator $\log_\theta\! B$, when
the spectral cut is at a general angle $\theta $. 
When $B$ has spectral cuts at two angles $\theta <\varphi $, one can
define the sectorial projection $\Pi _{\theta ,\varphi }(B)$ whose 
range contains the generalized eigenspaces for eigenvalues
with argument in $\,]\theta ,\varphi [\,$; this is studied in the
last part of the paper.
The operator $\Pi _{\theta ,\varphi }(B)$ is shown to be proportional 
to the difference between
$\log_\theta\! B$ and $\log_\varphi\! B$, having slightly better
symbol properties than they have. We show by examples that it 
belongs to the Boutet 
de Monvel calculus in many special cases, but lies
outside the calculus in general.
\endabstract
%\subjclassyear{2000}\subjclass 35S15, 58J42 \endsubjclass
\endtopmatter

\subhead 1. Introduction
\endsubhead

The purpose of this paper is to set up logarithms and sectorial projections
for elliptic boundary value problems, and to establish and analyze 
residue definitions associated with these operators. Let us
first recall the situation for boundaryless manifolds:

For a classical elliptic pseudodifferential operator ($\psi $do) $P$ of order
$m>0$, acting in a vector bundle $\widetilde E$ over a closed (i.e.,
compact boundaryless) $n$-dimensional manifold $\widetilde X$, certain
functions of
the operator have been studied with great interest for many
years. Assuming that $P$ has no
eigenvalues on some ray, say ${\Bbb R}_-$, one has from
Seeley's work \cite{S1} that the complex powers $P^{-s}$ can be
defined as $\psi $do's by use of the resolvent $(P-\lambda )^{-1}$. 
Moreover, the zeta function $\zeta
(P,s)=\Tr(P^{-s})$ has a meromorphic extension to $s\in{\Bbb C}$ with
at most simple poles at the real numbers $\{(n-j)/m\mid j\in{\Bbb
N}\}$ (we denote $\{0,1,2,\dots\}= {\Bbb
N}$). There is no pole at $s=0$ (for $j=n$), and the value $\zeta
(P,0)$ plays an important role in index formulas. Let us define the
{\it basic zeta value} $C_0(P)$ by 
$$
C_0(P)=\zeta (P,0)+\nu _0,\tag1.1
$$
where $\nu _0$ is the algebraic multiplicity of the zero eigenvalue of
$P$ (if any). It is well-known how $C_0(P)$ can be calculated 
in local coordinates from finitely many homogeneous terms of the
symbol of $P$. 

Another interesting function of $P$ is $\log P$, defined on smooth
functions by
$$
\log P=\lim _{s\searrow 0}\tfrac i{2\pi }\int_{\Cal C}\lambda
^{-s}\log \lambda \,(P-\lambda )^{-1} \,d\lambda;
\tag1.2
$$ 
here $\lambda ^{-s}$ and $\log \lambda $ are taken with branch cut
$\rmi$, and $\Cal C$
is a contour in ${\Bbb C}\setminus\crm$ going
around the nonzero spectrum of $P$ in the positive direction. By use
of the fact that $\log P=-\frac d{ds}P^{-s}|_{s=0}$, Scott \cite{Sc}
showed that 
$$
C_0(P)=-\tfrac1m \operatorname{res}(\log P),\tag1.3
$$
where $\operatorname{res}(\log P)$ is a slight generalization of 
Wodzicki's noncommutative residue (\cite{W2}, Guil\-le\-min  \cite{Gu}).

In the case of a compact $n$-dimensional manifold $X$ with boundary $\partial
X=X'$ (smoothly imbedded in an $n$-dimensional manifold
$\widetilde X$ without boundary), one can study the analogous
operators and constants
defined from a realization $B$ of a pseudodifferential (or
differential) elliptic boundary value problem. Here $B=(P+G)_T$,
defined from a system 
$\{P_++G ,T\}$ of order $m>0$ ($m\in{\Bbb Z}$) in the Boutet de Monvel calculus
\cite{B}, where $P$ is a $\psi $do on $\widetilde X$ and $P_+$ is its
truncation to $X$ (acting in $E=\widetilde E|_X$), $G$ is a singular Green operator (s.g.o.)\ and $T$
is a system of trace operators. $B$ is the operator
acting like
$P_++G$ with domain 
$$
D(B)=\{u\in H^m(X,E)\mid Tu=0\}, \tag1.4
$$
where $H^m(X,E)$ is the Sobolev space of order $m$.
In the differential operator case, $G=0$. Assuming that for $\lambda $ on a ray, say $\rmi$, $\{P_++G -\lambda
,T\}$ satisfies the hypotheses of parameter-ellipticity of Grubb
\cite{G1, Sect.\ 3.3} (consistent with those of Seeley \cite{S2} in
the differential operator case), one can define the complex powers by functional
analysis and study the pole structure of $\zeta (B,s)=\Tr(B^{-s})$ 
\cite{G1, Sect.\ 4.4}, and in particular discuss the basic zeta value
$C_0(B)$ defined similarly to (1.1).
However, in contrast with the closed manifold case, the powers
$B^{-s}$ do not lie in the calculus we are using (in particular their
$\psi $do part does not satisfy the transmission condition of
\cite{B}).
Then it is advantageous to build the analysis more directly on the
resolvent, which does belong to the parameter-dependent calculus set
up in \cite{G1}. In fact,  for $N>n/m$ (such that $(B-\lambda )^{-N}$
 is trace-class),
there is a trace expansion for $\lambda \to\infty $ 
in a sector $V$ around ${\Bbb R}_-$:
$$
\Tr(B-\lambda )^{-N}=\sum_{0\le j\le n}c^{(N)}_j(-\lambda
)^{(n-j)/m-N}
+O(\lambda ^{-N-\varepsilon })\tag1.5
$$
($\varepsilon >0$), and here
$$
C_0(B)=c^{(N)}_n ,\tag1.6
$$
independently of $N$. It is shown in \cite{G3} that for a
generalization of (1.3) to $B$,
$$
C_0(B)=-\tfrac1m\res(\log B),\tag1.7
$$
 it is sufficient to be able to define $\log B$; the complex powers 
$B^{-s}$ are not needed.

The present paper gives in Sections 2 and 3 a detailed study of $\log
B$. For one thing, this allows a more precise interpretation of the
formula (1.7),
initiated in \cite{G3}. Another important purpose is to 
open up for the use of compositions of
$\log B$ with other operators. These are needed for the
consideration of composed zeta functions $\zeta
(A,B,s)=\Tr(AB^{-s})$ with general $A$ from the calculus of \cite{B},
or rather, trace expansion formulas for composed resolvents 
$A(B-\lambda )^{-N}$. Such a study is
carried out in \cite{G4} using the results on $\log B$ obtained 
in the present paper. We show in Section 2 that
$$
\log B=(\log P)_+ + G^{\operatorname{log}},\tag1.8
$$
where $G^{\log}$ is a generalized singular Green operator satisfying a
specific part of the usual symbol estimates for s.g.o.s; its principal
part has a singularity at the boundary. In Section 3 we study its residue.

\medskip
If,
more generally than $\rmi$, the ray free of eigenvalues for $B$ (the
spectral cut) is $e^{i\theta }\rp$ for some angle $\theta $, 
the corresponding operator functions will be
defined by formulas where $\lambda ^{-s}$ and $\log \lambda $ (as in (1.2)) are
replaced by $\lambda _\theta ^{-s}$ and $\log _\theta \!\lambda $ with
branch cut $e^{i\theta }\rp$, and the integration curve runs in
${\Bbb C}\setminus e^{i\theta }\crp$. The functions 
are then provided with an index
$\theta $;
$$ 
\zeta_\theta (B,s)=\Tr(B_\theta ^{-s}),\quad
\log_\theta\! B =(\log _\theta \!P)_+ + G^{\operatorname{log}_\theta}.
\tag1.9
$$
When $B$ has spectral cuts at 
$\theta $ and $\varphi $ for some $\theta <\varphi <\theta +2\pi $, it
is of interest to study the {\it  sectorial
projection} $\Pi _{\theta ,\varphi }(B)$, a
projection whose range contains the generalized eigenspace of $B$ for the
sector
 $\Lambda _{\theta
,\varphi }=\{\, re^{i\omega }\mid r>0,\, \theta <\omega <\varphi \,\}$
and whose
nullspace contains the generalized eigenspace of $B$ for
 $\Lambda _{\varphi ,\theta +2\pi  }$; 
it was
considered earlier by Burak \cite{Bu}, and in the boundaryless case by
Wodzicki \cite{W2}, Ponge \cite{P}.
We show in Section 4 that it equals $\tfrac 
i{2\pi }(\log_\theta \!B -\log_\varphi \!B)$ and has the form 
$$
\Pi _{\theta ,\varphi }(B)= (\Pi _{\theta ,\varphi }(P))_+ 
+G_{\theta ,\varphi }.\tag1.10
$$
Here $\Pi _{\theta ,\varphi }(P)$ is a zero-order classical  $\psi $do,
which satisfies the transmission condition when $m$ is even,
and $G_{\theta ,\varphi }$ is a generalized s.g.o, bounded in $L_2$
in the differential operator case. There are natural types of examples
where $G_{\theta ,\varphi }$ is a standard s.g.o.\ as in
\cite{B}, but in general it will be of a generalized type satisfying only part of
the standard symbol estimates.
 
We expect to take up elsewhere the study of its residue, whose
possible vanishing is important for the study of eta functions
associated with $B$.

\subhead 2. The singular Green part of the logarithm
\endsubhead

Let $X$ be a compact $n$-dimensional $C^\infty $ manifold with
boundary $\partial X=X'$, provided with a hermitian $C ^\infty$ vector bundle
$E$. We can assume that  $X$ is smoothly imbedded in an
$n$-dimensional manifold $\widetilde X$ without boundary and that $E$ is
the restriction to $X$ of a bundle $\widetilde E$ over $\widetilde X$. 
Consider a system $\{P_++G,T\}$ of operators in the Boutet de Monvel 
calculus \cite{B}
(pseudodifferential boundary operators, $\psi $dbo's). Here $P$
is defined as a $\psi $do of order $m>0$ on 
$\widetilde X$ acting on the sections of $\widetilde E$,
and its truncation to $X$ is 
$$
P_+=r^+Pe^+,\quad r^+\text{ restricts from
$\widetilde X$ to $X$, }e^+\text{ extends by 0}.\tag2.1
$$
To assure that $P_+$ maps $C^\infty (X,E)$ into itself, $P$ is assumed
to satisfy the transmission condition, which means that
in local coordinate systems at the boundary, where the manifold is
replaced by $\rnp=\{x=(x_1,\dots,x_n)\mid x_n>0\}$, with notation
$x'=(x_1,\dots, x_{n-1})$,
$$
     \partial_x^\beta \partial_\xi ^\alpha p_{m-j}(x',0,0,-\xi _n)
=(-1)^{m-j-|\alpha |}\partial_x^\beta \partial_\xi ^\alpha
p_{m-j}(x',0,0,\xi _n)\text{ for }|\xi _n|\ge 1, \tag2.2
$$
for all indices; $m$ is integer. (A discussion of such conditions can
be found in Grubb and H\"o{}rmander \cite{GH}.) 
$G$ is a singular Green operator in $E$ of order and class $m$, 
and $T=\{T_0,\dots,T_{m-1}\}$ is a system of 
trace operators $T_k$ of order and class $k$, 
going from $E$ to
bundles $F_k$ over $\partial X$, defining an elliptic boundary
value problem.
In particular, 
$$\sum_{0\le k\le m-1}\dim F_k=\tfrac12
m\dim E.
\tag2.3$$
Details on these operator types can be found in  \cite{B}, \cite{G1}.

 We assume that the system $\{P_++G-\lambda ,T\}$
satisfies the conditions of
parameter-ellipticity in \cite{G1, Def.\ 3.3.1} for $\lambda
$ on the rays in a sector $V$ around  $\Bbb R_-$. 
In particular, it can be a differential operator system; here $P$ and
$T$ are differential, and $G$ is omitted. A
classical example is the Laplace operator on a domain in ${\Bbb R}^n$,
together
with the Dirichlet trace operator $T=\gamma _0$.
 
It should be noted that the hypotheses imply that the trace operator
is normal, as
accounted for in \cite{G1, Section 1.5}.

The system has a certain regularity number $\nu
$ in the sense of \cite{G1}; it is an integer or half-integer in
$[\frac12,m]$ for pseudodifferential problems, $+\infty $ for purely
differential problems. 

{}From the system we define the realization $B=(P+G)_T$ as the operator
acting like $P_++G$ with domain (1.4).
By \cite{G1, Ch. 3}, the resolvent
$R_\lambda =(B-\lambda )^{-1}$ exists on each ray in $V$ for
sufficiently large 
$|\lambda |$, and is $O(\lambda ^{-1})$ in $L_2$ operator norm there.
It has the structure$$
R_\lambda =Q_{\lambda ,+}+G_\lambda ,\tag2.4
$$
where $Q_\lambda =(P-\lambda )^{-1}$ on $\widetilde X$ (which can be
assumed to be compact), and $G_\lambda
$ is the singular Green part. Since the spectrum of $B$ is
discrete, we can assume (after a small rotation if necessary) that
$\Bbb R_-$ is free of eigenvalues of $B$, and likewise for $P$.
 
We shall define the operator $\log (B)=\log ((P+G)_T)$, also written
$\log B$, $\log
(P+G)_T$, by  
$$
\log (P+G)_T=\lim _{s\searrow 0}\tfrac i{2\pi }\int_{\Cal C}\lambda
^{-s}\log \lambda \,R_\lambda \,d\lambda,
\tag2.5
$$ 
to be further explained below;
here $\Cal C$ is a Laurent loop 
$$
\Cal C=\{re^{i\pi }\mid \infty >r>r_0\}\cup\{r_0e^{i\omega }\mid \pi \ge \omega \ge -\pi \}\cup\{re^{-i\pi
}\mid 
r_0<r<\infty \}\tag2.6
$$
going around the nonzero spectrum of $(P+G)_T$ in the positive direction.

Insertion of the decomposition (2.4) in the defining formula (2.5) shows
that $Q_{\lambda ,+}$ contributes with $$
\lim _{s\searrow 0}\tfrac i{2\pi }\int_{\Cal C}\lambda
^{-s}\log \lambda \,r^+Q_{\lambda}e^+ \,d\lambda=r^+(\log P) e^+=(\log P)_+,
\tag2.7
$$
where $\log P$ is well-known from the closed manifold case, cf.\ (1.2). Its symbol in local
coordinates is of the form $$
\operatorname{symb}(\log P)=m\log[\xi ]+l(x,\xi ),\tag2.8
$$
where $l(x,\xi )$ is a classical $\psi $do symbol of order 0 
(see also the lemma below), and 
$[\xi ]$ is a smooth positive function that equals $|\xi |$ for $|\xi
|\ge 1$. 
The operator is continuous from $H^t(\widetilde X,\widetilde E)$ to
$H^{t-\varepsilon }(\widetilde X,\widetilde E)$ for any $\varepsilon
>0$; hence$$
(\log P)_+\:H^t( X,E)\to H^{t-\varepsilon }( X,E)\text{ for
}t>-\tfrac12.\tag2.9 
$$
(The limit for $s\to 0$ in (2.7) can be taken in this operator norm.) 

In even-order cases, the transmission condition satisfied by $P$ carries over to
$l(x,\xi )$:

\proclaim{Lemma 2.1} When $m$ is even, $l(x,\xi )$ satisfies the
transmission condition.
\endproclaim

\demo{Proof} As shown e.g.\ in Okikiolu \cite{O}, the symbol of $\log
P$ is calculated in local coordinates from the symbol $q(x,\xi ,\lambda )$
of $Q_\lambda $ by integration with $\log \lambda $ around the
spectrum of the principal symbol $p_m$ of $P$; here the
quasi-homogeneous terms in the expansion $q(x,\xi ,\lambda )\sim
\sum_{j\in{\Bbb N}}q_{-m-j}(x,\xi ,\lambda )$ (homogeneous of degree
$-m-j$ in $(\xi ,|\lambda |^{\frac1m})$ on each ray) contribute as
follows:
$$\aligned
\tfrac i{2\pi }\int_{\Cal C(x,\xi )}&\log \lambda \,q_{-m}(x,\xi
,\lambda ) \,d\lambda=
\tfrac i{2\pi }\int_{\Cal C(x,\xi )}\log \lambda \,(p_m(x,\xi
)-\lambda )^{-1} \,d\lambda\\&=\log p_m(x,\xi )
=\log  ([\xi ]^m)+\log([\xi ]^{-m}p_m(x,\xi ))=m\log[\xi ]+l_0(x,\xi ),\\
\tfrac i{2\pi }\int_{\Cal C(x,\xi )}&\log \lambda \,q_{-m-j}(x,\xi
,\lambda ) \,d\lambda=
l_{-j}(x,\xi )\text{ for }j>0,\\
\endaligned\tag2.10$$
where $\Cal C(x,\xi )$ is a closed curve in $\Bbb C\setminus\crm$
around the spectrum of $p_m(x,\xi )$. Each $l_{-j}$ is homogeneous in
$\xi $ of degree $-j$ for $|\xi |\ge 1$; for $j=0$ it follows since
$[\xi ]^{-m}p_m(x,\xi )$ is so, and for $j\ge 1$ it is seen e.g.\ as
follows (where we set $\lambda =t^m\varrho $):
$$\aligned
l_{-j}(x,t\xi )&=\tfrac i{2\pi }\int_{\Cal C(x,t \xi )}\log \lambda
\,q_{-m-j}(x,t\xi ,\lambda ) \,d\lambda\\
&= \tfrac i{2\pi }\int_{t^{-m}\Cal C(x,t \xi )}(\log
\varrho+m\log t)t^{-m-j} \,q_{-m-j}(x,\xi
,\varrho  ) \,t^md\varrho \\ 
&= t^{-j}l_{-j}(x,\xi )+mt^{-j}\log t\tfrac i{2\pi }\int_{\Cal C(x,\xi
)} \,q_{-m-j}(x,\xi ,\varrho  ) \,d\varrho ,
\endaligned$$ 
where the last term is zero since $q_{-m-j}$ is $O(|\varrho |^{-2})$
for $|\varrho |\to\infty $ when $j>0$.

When $m$ is even, we see that the transmission condition (2.2) 
carries over through the calculations (2.10) to the corresponding
property for $l(x,\xi )$, since the parity of $-j$ is the same as that
of $-j-m$.\qed
\enddemo

Now consider the contribution from $G_\lambda $. Here we shall use the
following observations:
$$
\aligned
Q_\lambda +\lambda ^{-1}&=Q_\lambda +\lambda ^{-1}(P-\lambda
)Q_\lambda =\lambda ^{-1}PQ_\lambda \text{ on
}\widetilde X,\\ 
R_\lambda +\lambda ^{-1}&=R_\lambda +\lambda
^{-1}(P_++G-\lambda )R_\lambda\\
&=\lambda ^{-1}(P_++G)(Q_{\lambda ,+}+G_\lambda )\\
&=\lambda ^{-1}[(PQ_\lambda )_+-L(P,Q_\lambda )+GQ_{\lambda ,+}
+(P_+ +G)G_\lambda ]\\
&=Q_{\lambda,+} +\lambda ^{-1}+\lambda ^{-1}[-L(P,Q_\lambda
)+GQ_{\lambda ,+} +(P_+ +G)G_\lambda ] \text{ on }X;
\endaligned\tag2.11
$$ 
they imply in view of (2.4) that $G_\lambda $ may be written as $$
G_\lambda =\lambda ^{-1}[-L(P,Q_\lambda )+GQ_{\lambda ,+}
+(P_+ +G)G_\lambda ].\tag 2.12
$$
Here $L(P,Q_\lambda )=G^+(P)G^-(Q_\lambda )$ in local coordinates.
(The latter formula is accounted for in \cite{G1, (1.2.49--50) and
Sect.\ 2.6}; we recall that $G^+( P)=r^+ Pe^-J$
and $G^-( P)=Jr^- Pe^+$, where $e^\pm$ extends by zero from
$\rnpm$ to ${\Bbb R}^n$, $r^\pm$ restricts from ${\Bbb R}^n$ to
$\rnpm$, and $J$ is the reflection map $J\: u(x',x_n)\mapsto u(x',-x_n)$.) 
By \cite{G1, Th.\ 3.3.2}, $G_\lambda $ is of order $-m$ and
regularity $\nu $; moreover, (2.12) shows that it is $\lambda ^{-1}$
times an s.g.o.\ of order 0 and regularity $\nu $ (by the composition
rules in \cite{G1, Th.\ 2.7.6--7}).

Since $$
Q_\lambda \:L_2(\widetilde X,\widetilde E)\to H^{m-\varepsilon
}(\widetilde X,\widetilde E),\quad 
G_\lambda \:L_2( X, E)\to H^{m-\varepsilon
}( X, E),\text{ with norms }O(\lambda ^{-\varepsilon /m}),
$$
for $\varepsilon\in [0,m]$ (a standard observation used also in \cite{G1,
pp.\ 409--410}), each of the terms in [\;] in (2.12) maps
$L_2(X,E)$ to $H^{-\varepsilon }(X,E)$ with norm $O(\lambda
^{-\varepsilon /m})$. Then we can perform the integration in this
operator norm (letting $s\to 0)$, defining the s.g.o.-like part
$G^{\operatorname{log}}$ of $\log(P+G)_T$ by
$$
G^{\operatorname{log}}=\tfrac i{2\pi }\int_{\Cal C}\log \lambda
\,G_\lambda \,d\lambda
=\tfrac i{2\pi }\int_{\Cal C}\lambda ^{-1}\log \lambda
\,[-L(P,Q_\lambda )+GQ_{\lambda ,+}
+(P_+ +G)G_\lambda ] \,d\lambda,
\tag2.13
$$
also written as$$
\multline
G^{\operatorname{log}}=-G^+(P)\tfrac i{2\pi }\int_{\Cal C}\lambda ^{-1}\log \lambda
\,G^-(Q_\lambda )\,d\lambda \\
+G\tfrac i{2\pi }\int_{\Cal C}\lambda ^{-1}\log \lambda
\,Q_{\lambda,+}\,d\lambda 
+(P_++G)\tfrac i{2\pi }\int_{\Cal C}\lambda ^{-1}\log \lambda
\,G_{\lambda}\,d\lambda ,
\endmultline\tag2.14
$$
when localized.  
It is a bounded operator from 
$L_2( X, E)$ to $H^{-\varepsilon
}( X, E)$. Summing up, we have found:

\proclaim{Theorem 2.2} 
The logarithm of the realization $B=(P+G)_T$ satisfies
$$
\log B=\log(P+G)_T=(\log P)_+ + G^{\operatorname{log}},\tag2.15
$$
where $\log P$ is the logarithm of $P$ on $\widetilde X$, and
$G^{\log} $ is defined by {\rm (2.13), (2.14)}; the terms are 
bounded operators from 
$L_2( X, E)$ to $H^{-\varepsilon
}( X, E)$ (any $\varepsilon >0$). 
\endproclaim

The operator $G^{\operatorname{log}}$ is a generalized singular Green
operator, in the same spirit as 
the generalized s.g.o.s $G^{(-s)}$ studied in \cite{G1, Sect.\ 4.4} 
(the s.g.o.-like parts of the powers $B^{-s}$), and one 
can show as in \cite{G1, Th.\ 4.4.4} that there is a symbol-kernel
satisfying part of the usual $L_{2,x_n,y_n}(\Bbb R^2_{++})$ estimates
for s.g.o.s, allowing $D _{x'}^\beta $, $D _{\xi '}^\alpha
$, $(x_nD _{x_n})^k$ and $(y_nD _{y_n})^l$ in
arbitrarily high
powers (with exceptions for the principal term), and allowing some applications of $x_n^kD _{x_n}^{k'}$ and
$y_n^lD _{y_n}^{l'}$, limited by the regularity and other
restrictions. We account for this in Theorem 2.6 below; let us first
consider an example.

\example {Example 2.3} Let $P=1-\Delta $ on $\rnp$. It is easy to see
that the solution operator for the Dirichlet problem for $P-\lambda
=1-\Delta -\lambda $, $\lambda \in V={\Bbb C}\setminus\rp$, is
$R_\lambda =Q_{\lambda ,+}+G_\lambda $, where $Q_\lambda $ is the $\psi $do
$(1-\lambda -\Delta )^{-1}$ with symbol $(\ang\xi^2-\lambda )^{-1}$,
and $G_\lambda $ is the
singular Green 
operator with symbol-kernel $\frac{-1}{2\kappa _1}e^{-\kappa _1(x_n+y_n)}$;
$\kappa _1=(\ang{\xi '}^2-\lambda )^{\frac12}$. (We here use the
well-known notation $\ang{x}=(x_1^2+\dots+x_n^2+1)^{\frac12}$.) 
It follows that $$
\log P=\OP(2\log\ang\xi ).\tag2.16
$$
To find out how $G^{\operatorname{log}}$ acts on functions $\varphi
\in C_0^\infty (\rnp)$, we write (using that
$e^{-\kappa _1(x_n+y_n)}$ is rapidly decreasing in $\lambda $ on the
rays in $V$ when $y_n$
is in the support of $\varphi $):
$$\aligned
G^{\log}\varphi &=\tfrac i{2\pi }\int_{\Cal C}\log \lambda
\,G_{\lambda}\varphi \,d\lambda \\
&=\tfrac i{2\pi }\int_{\Cal C}\int_{{\Bbb R}^{n-1}}\int_0^\infty \log \lambda
\,e^{ix'\cdot \xi '}\tfrac{-1}{2\kappa _1}e^{-\kappa _1(x_n+y_n)}\acute \varphi (\xi ',y_n)\,dy_n\d \xi 'd\lambda ,
\endaligned$$
with $\acute \varphi $ denoting the partial Fourier transform $\acute
\varphi (\xi ',y_n)=\Cal F_{y'\to \xi '}\varphi(y',y_n)$. Here we can
calculate
$$  
\aligned
\tfrac i{2\pi }\int_{\Cal C}&\log\lambda \,
\tfrac
{-1}{2\kappa _1
}e^{-\kappa _1(x_n+y_n) }\,d\lambda 
=\int_{-\infty }^0\frac1{2(\ang{\xi '}^2-t)^{\frac12}} e^{-(\ang{\xi
'}^2-t)^{\frac12}(x_n+y_n)}\, dt\\
&=\int_{0 }^\infty \frac1{2(\ang{\xi '}^2+s)^{\frac12}} e^{-(\ang{\xi
'}^2+s)^{\frac12}(x_n+y_n)}\, ds
=\int_{\ang{\xi '} }^\infty \frac1{2u} e^{-u(x_n+y_n)}\, 2u\,du\\
&=\frac1{x_n+y_n}e^{-\ang{\xi '}(x_n+y_n)},
\endaligned\tag2.17$$ 
using that the $\log|\lambda |$ contributions cancel out (as
in \cite{G3, Lemma 1.2}).  Thus
$$
G^{\log} \varphi =\int_{{\Bbb R}^{n-1}}\int_0^\infty e^{ix'\cdot \xi
'}\frac1{x_n+y_n }e^{-\ang{\xi '} (x_n+y_n)}\acute \varphi (\xi
',y_n)\,dy_n\d \xi '
.
$$
This shows that $G^{\log}$ is a generalized kind of s.g.o.\ with symbol-kernel
$$
\tilde g^{\operatorname{log}}(x',x_n,y_n,\xi ')=\frac1{x_n+y_n}e^{-\ang{\xi
'}(x_n+y_n)}.\tag2.18
$$
Since the operator with kernel $\frac1{x_n+y_n}$ is bounded in
$L_2({\Bbb R}_+)$ (as a truncation of the Hilbert transform), it follows
that $G^{\log}$ is a bounded operator in $L_2(\rnp)$. 

Note that $\partial _{\xi _1}\tilde g^{\operatorname{log}}$ is a
standard s.g.o.\ symbol-kernel, and that $x_n\tilde
g^{\operatorname{log}}$ is  bounded.

The same calculations with $\ang{\xi '}$ replaced by $|\xi '|$ show
that for $P=-\Delta $, $G^{\log}$ has symbol-kernel $\frac1{x_n+y_n}
e^{-|\xi
'|(x_n+y_n)}$ for $|\xi '|\ge 1$.
\endexample

In the general {\it differential operator} case, $G^{\log}$ is qualitatively very much
like in this example. Here one can directly use the symbol-kernel estimates
and boundedness considerations worked out by Seeley in \cite{S2},
\cite{S3}. Notationally, we follow \cite{G3}; in particular,
the enumeration of quasi-homogeneous (resp.\ homogeneous) terms in the asymptotic
expansions of singular
Green symbol-kernels (resp.\ symbols) have been shifted by one step in comparison with \cite{G1}, in order to have
the same index on an s.g.o.\ symbol-kernel (resp.\ symbol) and its normal 
trace. For example, the principal part of a symbol-kernel $\tilde g$ of order
$-m$ is denoted $\tilde g_{-m}$ (although the corresponding symbol
$g_{-m}$ has
homogeneity degree $-m-1$). We shall use the notation $\dot\le$
(resp.\ $\dot\ge$) to
indicate ``less than or equal (resp.\ greater than or equal) to a
constant times'', and $\dot=$ to indicate that both $\dot\le$ and
$\dot\ge$ hold.

\proclaim{Theorem 2.4} Consider the case where $P$ is a differential
operator, $G=0$, and the trace operators $T_0,\dots, T _{m-1}$ are
differential operators. In this case, the singular Green part
$G_\lambda $ of the resolvent is of regularity $+\infty $ and its
symbol-kernel in local coordinates $\tilde g\sim \sum_{j\ge 0}\tilde g_{-m-j}$, expanded
in quasi-homogeneous terms
$$
\tilde
g_{-m-j}(x',\tfrac{x_n}t,\tfrac{y_n}t,t\xi',t^m \lambda )=t^{-m+1-j}\tilde g_{-m-j}(x',x_n,y_n,\xi', \lambda )
\text{ for }t\geq1,|\xi'|\geq1,\tag2.19
$$
 satisfies estimates on the rays in $V$, with $ \kappa =| \xi '|+|\lambda |^{\frac1m}$:
$$
|D_{x'}^\beta D_{\xi'}^\alpha x_n^kD_{x_n}^{k'}
y_n^lD_{y_n}^{l'}D_\lambda ^p
\tilde g_{-m-j}|\leg \kappa ^{1-m-|\alpha |-k+k'-l+l'-j-mp}e^{-c\kappa (x_n+y_x)}\tag2.20
$$
for all indices, when $ \kappa \ge \varepsilon $.

 Then $G^{\operatorname{log}}$ is, in local coordinates near $X'$, a
generalized singular Green operator
 $$
\aligned
G^{\operatorname{log}}u(x)
&=\int_{\Bbb R^{n-1}}\int_0^\infty
e^{ix'\cdot\xi'}\tilde g^{\operatorname{log}}(x',x_n,y_n,\xi')\acute
u(\xi',y_n)\,dy_n \d\xi'\\
&=\operatorname{OPG}(\tilde g^{\operatorname{log}}(x',x_n,y_n,\xi')) u(x)
\endaligned
\tag2.21
$$
with $\tilde g^{\operatorname{log}}\sim\sum_{j\in\Bbb N}\tilde g_{-j}^{\operatorname{log}}$;
here the $j$'th term 
is quasihomogeneous:
 $$
\tilde
g_{-j}^{\operatorname{log}}\left(x',\tfrac{x_n}t,\tfrac{y_n}t,t\xi'\right)
=t^{1-j}\tilde g_{-j}^{\operatorname{log}}(x',x_n,y_n,\xi')
\text{ for }t\geq1
\text{ and }|\xi'|\geq1,
\tag2.22
$$
and satisfies, when $|\xi '|\ge \varepsilon $,
$$|D_{x'}^\beta D_{\xi'}^\alpha
x_n^kD_{x_n}^{k'}y_n^lD_{y_n}^{l'}\tilde
g_{-j}^{\operatorname{log}}|
 \leg|\xi'|^{-|\alpha|-k+k'-l+l'-j}\tfrac1{x_n+y_n}e^{-c|\xi '|(x_n+y_n)}\tag2.23  
$$
for the indices satisfying
 $$
-k+k'-l+l'-|\alpha|-j\le 0.\tag2.24
$$

It follows in particular that $G^{\log}$ is a bounded operator in $L_p(X,E)$ for $1<p<\infty $.

\endproclaim 

\demo{Proof} The estimates (2.20) were shown in \cite{S2, (29)}, \cite{S3}. 
Because of the fall-off in $\lambda $, they allow us to define the
$j$'th term in the symbol-kernel of $G^{\log}$ for $|\xi '|\ge
\varepsilon $ by
 $$
\aligned\tilde g^{\operatorname{log}}_{-j}(x',x_n,y_n,\xi')&=\tfrac
i{2\pi}\int_{\Cal C}\log \lambda\,\tilde
g_{-m-j}(x',x_n,y_n,\xi',\lambda)\,d\lambda\\
&=\int_0^\infty \tilde
g_{-m-j}(x',x_n,y_n,\xi',-s)\,ds;
\endaligned\tag2.25
$$
here we rewrote the  
integral as in (2.17) (and \cite{G3, Lemma 1.2}). The homogeneity
is seen from the last integral, using (2.19). The function is estimated as follows, for the indices
satisfying (2.24), when we use that
$|\xi '|+s^{\frac1m}\eg (|\xi '|^m+s)^{\frac1m}$: 
  $$
\aligned
&|D_{x'}^\beta D_{\xi'}^\alpha
x_n^kD_{x_n}^{k'}y_n^lD_{y_n}^{l'}\tilde
g_{-j}^{\operatorname{log}}|=|\int_0^\infty D_{x'}^\beta D_{\xi'}^\alpha
x_n^kD_{x_n}^{k'}y_n^lD_{y_n}^{l'} \tilde
g_{-m-j}(x',x_n,y_n,\xi',-s)\,ds|\\
&\leg |\xi '|^{-|\alpha |-k+k'-l+l'-j}\int_0^\infty ((|\xi
'|^m+s)^{\frac1m})^{1-m}\
e^{-c (|\xi
'|^m+s)^{\frac1m}(x_n+y_n)}\, ds\\
&=
|\xi '|^{-|\alpha |-k+k'-l+l'-j}\int_{|\xi '|}^\infty  u ^{1-m}
e^{-c u(x_n+y_n)}mu^{m-1}\, du\\
& =|\xi'|^{-|\alpha|-k+k'-l+l'-j}\tfrac m{c(x_n+y_n)}e^{-c|\xi '|(x_n+y_n)}.\endaligned\tag2.26  
$$

The operator $G^{\log}$ is defined from a finite number of these
symbol terms multiplied with an excision function $ \zeta (|\xi '|)$,
where
$$
   \zeta (t)\in C^\infty ({\Bbb R}),\;\zeta (t)=0\text{ for }|t|\le
\delta _1,\;\zeta (t)=1\text{ for }|t|\ge
\delta _2,\tag 2.27
$$
plus an integral as in (2.13) of the remainder of
$G_\lambda $, which can be taken with arbitrarily high smoothness of
the kernel and
decrease for $\lambda \to\infty $, cf.\ \cite{S3, (2.14)}.
Applying the arguments of Theorem 1 of \cite{S3} (using Lemmas 1 and
2 there invoking Mihlin's theorem and the Hilbert transform) one finds 
that $G^{\log}$ is $L_p$-continuous as asserted.
\qed
\enddemo

\example{Remark 2.5} The lower order terms in $\tilde g^{\log}$ and
the derivatives are
not as singular for $x_n+y_n\to 0$ as (2.23) indicates. In fact,
the symbol-kernels one step down can be estimated as follows:
 $$
\aligned
\text{ When }-k+k'&-l+l'-|\alpha|-j\le -1,\\
|D_{x'}^\beta D_{\xi'}^\alpha
x_n^kD_{x_n}^{k'}y_n^lD_{y_n}^{l'}\tilde
g_{-j}^{\operatorname{log}}|
&\leg |\xi '|^{-|\alpha |-k+k'-l+l'-j+1 }\int_{|\xi
'|}^\infty 
u^{-1-\varepsilon }u^\varepsilon e^{-cu(x_n+y_n)}\, du\\
& \leg|\xi'|^{-|\alpha|-k+k'-l+l'-j+1+\varepsilon }\sup_{u\in {\Bbb
R}_+}|u^\varepsilon e^{-cu(x_n+y_n)}|\\
&\leg |\xi'|^{-|\alpha|-k+k'-l+l'-j+1+\varepsilon }
(x_n+y_n)^{-\varepsilon },
\endaligned\tag2.28  
$$
for $\varepsilon >0$. The symbol-kernels two steps down are bounded for $x_n+y_n\to
0$:
 $$
\aligned
\text{ When }-k+k'&-l+l'-|\alpha|-j\le -2,\\
|D_{x'}^\beta D_{\xi'}^\alpha
x_n^kD_{x_n}^{k'}y_n^lD_{y_n}^{l'}\tilde
g_{-j}^{\operatorname{log}}|
&\leg |\xi '|^{-|\alpha |-k+k'-l+l'-j+2}\int_0^\infty 
(|\xi'|+s^{\frac1m})^{-m-1}
\, ds\\
& \leg|\xi'|^{-|\alpha|-k+k'-l+l'-j+1},
\endaligned\tag2.29  
$$
and the smoothness at 0 increases with increasing $|\alpha |$ and $j$.
\endexample

Now let us turn to the pseudodifferential case and the methods of
\cite{G1, Sect.\ 4.4}.

\proclaim{Theorem 2.6} Let $\{P_++G,T\}$ have regularity
$\nu\in [\tfrac12,\infty[\,$,  and 
define $G^{\operatorname{log}}$ by {\rm (2.13)}.
 Then $G^{\operatorname{log}}$ is, in local coordinates near $X'$, a
generalized singular Green operator as in {\rm (2.21)}
with $\tilde g^{\operatorname{log}}\sim\sum_{j\in\Bbb N}
\tilde g_{-j}^{\operatorname{log}}$;
here the $j$'th term is quasihomogeneous as in  {\rm (2.22)} when
$j>0$,
and the series approximates $\tilde g^{\operatorname{log}}$
asymptotically in the sense that
 $$\|D_{x'}^\beta D_{\xi'}^\alpha
x_n^kD_{x_n}^{k'}y_n^lD_{y_n}^{l'}[\tilde g^{\operatorname{log}}
-\sum_{j<J}\tilde
g_{-j}^{\operatorname{log}}]\|_{L_{2,x_n,y_n}}
 \leg\ang{\xi'}^{-|\alpha|-k+k'-l+l'-J}\tag2.30  
$$
holds for the indices satisfying
 $$\aligned
-k+k'-l+l'-|\alpha|-J&<0,\\
[k-k']_-+[l-l']_-&<\nu.\endaligned\tag2.31
$$
Moreover, 
$$\|D_{x'}^\beta D_{\xi'}^\alpha
x_n^kD_{x_n}^{k'}y_n^lD_{y_n}^{l'}\tilde
g_{-J}^{\operatorname{log}}\|_{L_{2,x_n,y_n}}
 \leg\ang{\xi'}^{-|\alpha|-k+k'-l+l'-J}\tag2.32  
$$
holds for these indices.

With $\zeta (t)$ defined as in {\rm (2.27)}, the above
symbol-kernels multiplied with $\zeta (x_n)\zeta (y_n)$ satisfy
estimates for all $\alpha ,\beta ,J,k,k',l,l'$ with $\ang{\xi '}^{-M}$,
any $M$, in the right-hand side.
\endproclaim

\demo{Proof} This is modeled after the proof of \cite{G1, Th.\ 4.4.4}
and the remarks preceding it.

We recall from \cite{G1, Th.\ 3.3.9} that the symbol-kernel $\tilde
g(x',x_n,y_n,\xi ',\lambda )$ of $G_\lambda  $ (in a local coordinate
system) has an expansion in quasi-homogeneous terms $\tilde g\sim\sum
_{j\ge 0}\tilde g_{-m-j}$ satisfying (2.19) in $V$,
and that one has for all indices, denoting $\lambda =-\mu ^me^{i\omega
}$ ($\mu >0$), 
$(|\xi '|^2+\mu ^2+1)^{\frac12}=\ang{\xi ',\mu }$:
 $$\aligned
\|D_{x',\omega }^\beta&D_{\xi'}^\alpha x_n^kD_{x_n}^{k'}
y_n^lD_{y_n}^{l'}[\tilde
g-\sum_{j<J}\tilde g_{-m-j}]\|_{L_{2,x_n,y_n}}\\
&\leg
(\ang{\xi'}^{\nu-M'}+\ang{\xi',\mu}^{\nu-M'})
\ang{\xi',\mu}^{-m-\nu+M''}\\
&\leg\cases
\ang{\xi',\mu}^{-m-M'+M''},&\text{ when }M'\leq\nu,\\
\ang{\xi'}^{\nu-M'}\ang{\xi',\mu}^{-m-\nu+M''}&\text{ when
}M'\geq\nu,\endcases\endaligned\tag2.33 
$$
with
 $$\aligned
M'&=[k-k']_++[l-l']_++|\alpha|+J,\\
M''&=[k-k']_-+[l-l']_-\;;\;\text{ so
}\\
-M'+M''&=-k+k'-l+l'-|\alpha|-J.\endaligned\tag2.34
$$
The notation $N_{\pm}=\max\{\pm N,0\}$ is used, and we have (as
recalled earlier) changed the
indexation from \cite{G1} by one step as in \cite{G3}.

Let us first observe that the \lq\lq error terms" and remainders in
the resolvent 
construction, that are negligible in the class of operators of order
$-m$ and regularity $\nu$, give rise to generalized s.g.o\. error
terms $G'$ here, satisfying estimates of the type (as in \cite{G1,
Lemma 2.3.11})
$$\multline
\|D_{x'}^\beta D_{\xi'}^\alpha  x_n^kD_{x_n}^{k'}y_n^lD_{y_n}^{l'}\tilde
g'\|_{L_{2,x_n,y_n}} \leg
\ang{\xi'}^{-M}|\int_{\Cal C}\log\lambda\,
\ang\lambda^{-1-(\nu-[k-k']_--[l-l']_-)/m}\,d\lambda|\\
 \leg\ang{\xi'}^{-M},\text{ for any }M, \text{
when }[k-k']_-+[l-l']_-< \nu.\endmultline
\tag2.35 
$$
It follows that the corresponding kernels $\Cal K_{G'}(x,y)$
satisfy, for these indices:
$$\sup_{x',y'}\|D_{x',y'}^\gamma x_n^kD_{x_n}^{k'}y_n^lD_{y_n}^{l'}\Cal
K_{G'}\|_{L_{2,x_n,y_n}}<\infty.\tag2.36
$$ 

For $j>0$ the $L_{2,x_n,y_n}$-norm of $\tilde g_{-m-j}$ is $O(\lambda
^{-1-1/2m})$ since $\nu \ge \frac12$, so the corresponding term $\tilde
g_{-j}^{\operatorname{log}}$ can be defined directly for $|\xi '|\ge
1$ by Cauchy integrals as in (2.25), convergent in the
$L_{2,x_n,y_n}$-norm. 
The quasi-homogeneity of $\tilde g_{-j}^{\operatorname{log}}$ is seen
as in (2.25) by using \cite{G3, Lemma 1.2} in $L_{2,x_n,y_n}$-norm.

% The terms in $\tilde g$ with $j>0$ are $O(\lambda ^{-1-1/2m})$ in
% $L_{2,x_n,y_n}$-norm since
% $\nu \ge \frac12$, so the corresponding terms in $\tilde
% g^{\operatorname{log}}$ can be defined directly for $|\xi '|\ge 1$
% by Cauchy integrals as in (2.25), convergent in the
% $L_{2,x_n,y_n}$-norm. 
% The last identity in (2.25) follows by using
% \cite{G3, Lemma 1.2} in $L_{2,x_n,y_n}$-norm, and the homogeneity is
% seen from this.

We use the estimates (2.33) to see that for $\tilde
g^{\operatorname{log}}-\sum_{j<J}\tilde g^{\operatorname{log}}_{-j}$
with $J>0$ (so that the first term is excluded),  the integrand in
the corresponding 
Cauchy integral is $O(\lambda^{-1-\varepsilon})$ in
$L_{2,x_n,y_n}$-norm (some $\varepsilon >0$), when
 $$
-k+k'-l+l'-|\alpha|-J<0,
\text{ if
}[k-k']_++[l-l']_++|\alpha|+J\leq\nu,
\tag2.37
$$
and when
 $$
[k-k']_-+[l-l']_-<\nu,
\text{ if
}[k-k']_++[l-l']_++|\alpha|+J\geq\nu.
\tag2.38
$$
Then the integral converges and defines a symbol-kernel satisfying
the asserted estimate. Since
 $$
-k+k'-l+l'-|\alpha|-J
=[k-k']_-+[l-l']_--([k-k']_++[l-l']_++|\alpha|+J),
 $$
 we see that the conditions \lq\lq if $\dots$" can be left out in
(2.37)--(2.38), leading to the formulation (2.31).

We still have to consider the first term $\tilde g^{\operatorname{log}}_0$
in $\tilde g^{\operatorname{log}}$,
defined from the principal part $\tilde g_{-m}$ of $\tilde g$. Here we
use that 
$\tilde g_{-m}$ can be found by performing the resolvent construction
on the principal boundary symbol level for
the corresponding
operators on $L_2(\rp)$, and that they obey a one-dimensional version of
the identities in (2.11). So we can replace $\tilde g_{-m}$ by the
symbol-kernel of the principal boundary symbol version of (2.12), which gives a
convergent Cauchy integral, when the $\lambda $-independent factors
are pulled outside of the integration. In a formal sense, we can ascribe it a
symbol-kernel $\tilde g^{\operatorname{log}}_0(x',x_n,y_n,\xi ')$. 
The resulting boundary 
symbol operator is  
continuous from $L_2(\rp)$ to $H^{-\varepsilon }(\rp)$ for
$\varepsilon >0$, at each $(x',\xi ')$.
If we define the functions derived from $\tilde g^{\operatorname{log}}_0$
``weakly''
by
$$  \multline
D_{x'}^\beta D_{\xi'}^\alpha
x_n^kD_{x_n}^{k'}y_n^lD_{y_n}^{l'}\tilde
g_{0}^{\operatorname{log}}(x',x_n,y_n,\xi ')\\
=\tfrac
i{2\pi}\int_{\Cal C}\log \lambda\,D_{x'}^\beta D_{\xi'}^\alpha
x_n^kD_{x_n}^{k'}y_n^lD_{y_n}^{l'}\tilde
g_{-m}(x',x_n,y_n,\xi',\lambda)\,d\lambda,\endmultline
$$ 
we can use that the integral converges in $L_{2,x_n,y_n}$-norm when
the indices 
satisfy (2.31). In this sense, the estimates (2.30) hold also when
$J=0$ in (2.31). 

The estimates (2.32) of the individual terms follow from
(2.30) since $ \tilde
g_{-J}^{\operatorname{log}}=$ \linebreak $(\tilde g^{\operatorname{log}}-\sum_{j<J}\tilde
g_{-j}^{\operatorname{log}})- (\tilde g^{\operatorname{log}}-\sum_{j<J+1}\tilde
g_{-j}^{\operatorname{log}})$. 

Finally, for the statements on the symbol-kernels multiplied with
$\zeta (x_n)\zeta (y_n)$, note that $\zeta (t)$ can for any $k\in\Bbb
N$ be written as
$t^k\zeta _k(t)$ with a bounded smooth function $\zeta _k$, so from the
already shown estimates we can infer arbitrarily rapid fall-off in
$\xi '$ by
rewriting with arbitrarily high powers of $x_n$ and $y_n$.  
\qed\enddemo

If $R_\lambda $ has infinite regularity, $\nu $ can be arbitrarily large
in the second line of (2.31), so the line can be left out. Note that even
then there is a limitation on the indices for which we get standard
s.g.o.\ estimates. 

\medskip
While $G^{\log}$ is the primary s.g.o.-type operator to consider in
this connection, it is also of interest to study some other s.g.o.-type
operators here, namely, in local coordinates, $G^+(\log P)=r^+(\log P)e^-J$
and $G^-(\log P)=Jr^-(\log P)e^+$, with notation as in the text after (2.12).
The operators $G^\pm (\log P)$ have properties very similar to those
of $G^{\operatorname{log}}$:

\proclaim{Theorem 2.7} The operators $G^\pm(\log P)$ are defined in local
coordinates by$$
\aligned
G^+(\log P)&=r^+\log Pe^-J=r^+\tfrac i{2\pi }\int_{\Cal C}\log \lambda
\,Q_\lambda \,d\lambda\, e^-J\\
&=\tfrac i{2\pi }\int_{\Cal C}\lambda ^{-1}\log \lambda
\,G^+(PQ_{\lambda }) \,d\lambda,\\
G^-(\log P)&=Jr^-\log Pe^+=Jr^-\tfrac i{2\pi }\int_{\Cal C}\log \lambda
\,Q_\lambda \,d\lambda\, e^+\\
&=\tfrac i{2\pi }\int_{\Cal C}\lambda ^{-1}\log \lambda
\,G^-(PQ_{\lambda })\,d\lambda
.\endaligned\tag2.39
$$
Their symbol-kernels $\tilde g^{\pm}(\operatorname{log}p)$ have
properties like those 
of $\tilde g^{\operatorname{log}}$ in Theorem {\rm 2.6}, with $\nu
=m$.

In particular, when $P$ is a differential operator, the s.g.o.s
$G^{\pm}(Q_\lambda )$ satisfy Seeley's estimates {\rm (2.20)}, and
hence the operators $G^{\pm}(\log P)$ have symbol estimates and
boundedness properties like those
of $G^{\log}$ in Theorem {\rm 2.4}, Remark {\rm 2.5}.
\endproclaim 

\demo{Proof} The defining integrals are established by use of the
first formula in (2.11), noting that $G^{\pm}(\lambda ^{-1})=0$. 
By \cite{G1, Th.\ 2.7.4}, $G^\pm(Q_\lambda )$ is a parameter-dependent
polyhomogeneous family of s.g.o.s of order $-m$ and regularity
$m-\varepsilon $ (any $\varepsilon >0$), since $Q_\lambda $ is of
order $-m$ and regularity $m$. The symbol-kernel then satisfies
estimates like those for $\tilde g$ in Theorem 2.6, with $\nu
=m-\varepsilon $. The method of Theorem 2.6 leads to the conclusion
that the resulting symbol-kernel $\tilde g^\pm(\log p)$ has properties like
those stated for $\tilde g^{\log}$, with $\nu =m-\varepsilon $; here
$\varepsilon $ can be removed since the second inequality in (2.31) is
sharp. 

For the second statement, we must show that the Seeley estimates
(2.20) are valid for the homogeneous terms in the symbol-kernel of
$G^\pm (Q_\lambda )$. But this is easy. Consider e.g.\ $G^+(Q_\lambda
)$. Using the Taylor
expansion of the symbol of $Q_\lambda $ at $x_n=0$:
$$
q(x',x_n,\xi ,\lambda )\sim \sum _{l\in{\Bbb
N}}\tfrac1{l!}x_n^l\partial_{x_n}^lq(x',0,\xi ,\lambda )
$$
we have from \cite{G1, Th.\ 2.7.4} that
$$
g^+(q)(x',\xi ,\eta _n,\lambda )\sim \sum_{l\in{\Bbb
N}}\tfrac1{l!}\overline D_{\xi _n}^l g^+[\partial_{x_n}^lq(x',0,\xi
,\lambda )],  
$$
where $g^+[f](\xi _n,\eta _n)$ is the s.g.o.\ symbol corresponding to the
symbol-kernel $\tilde g^+[f](x_n,y_n)$ defined by:
$$
\tilde g^+[f](x_n,y_n)=\bigl(r^+_{z_n}[\Cal F_{\xi _n\to z_n}^{-1}f]\bigr)|_{z_n=x_n+y_n}.
$$
The homogeneous terms in the symbols $\partial_{x_n}^lq(x',0,\xi
,\lambda )$
are rational functions of $\xi _n$ with $\frac12 m \dim E$ poles in 
${\Bbb C}_\pm=\{z\in{\Bbb C}\mid \operatorname{Im}z\gtrless 0\}$, lying inside a circle of radius $C\kappa $ and having
a distance $\ge c\kappa $ from the real axis, for suitable positive
constants $C>c$. (A more detailed description is given e.g.\ in
\cite{G1, Remark 3.3.7}.) For simplicity of notation, consider the $j$'th
term $q_{-m-j}$ itself.
The inverse Fourier transform evaluated at
$z_n>0$ can be written as an integral of $e^{iz_n\xi
_n}q_{-m-j}(x',0,\xi ',\xi _n)$ over the curve bounding
 the intersection of the
circle $\{|\xi _n|=C\kappa \}$ with the halfplane $\{\operatorname{Im}\xi _n\ge
c\kappa \}$ (lying in ${\Bbb C}_+$). 
We get the factor $e^{-c\kappa z_n}$ since $|e^{iz_n\xi _n}|\le
e^{-c\kappa z_n}$ on the curve. (Similarly, the inverse Fourier
transform evaluated at $z_n<0$ can be written as an integral over a
closed curve in ${\Bbb C}_-$ with $\operatorname{Im}\xi _n\le -c\kappa
$.) For the resulting symbol-kernel, this gives the factor
$e^{-c\kappa (x_n+y_n)}$; the power of $\kappa $ in front is seen
from the degree of the rational function.

Once the estimates (2.20) are established, the rest of the proof goes
as in Theorem 2.4. 
\qed
\enddemo 

\example{Example 2.8}
For $P=1-\Delta $ as in Example 2.3, one finds by direct calculation
of the inverse Fourier transform w.r.t.\ $\xi _n$ that 
 $G^\pm(Q_\lambda )$ both
have the symbol-kernel$$
\tilde g^+=\tilde g^-=\tfrac{1}{2\kappa _1}e^{-\kappa _1(x_n+y_n)},\tag 2.40
$$
with $\kappa _1=(\ang{\xi '}^2-\lambda )^{\frac12}$. 
Then the calculations of Example 2.3 can be used again, to see that  
$$
\tilde g^{+}(\operatorname{log}p)
(x',x_n,y_n,\xi ')
=\tilde g^{-}(\operatorname{log}p)
(x',x_n,y_n,\xi ')=\frac{-1}{x_n+y_n}e^{-\ang{\xi
'}(x_n+y_n)}
.\tag 2.41
$$

For $P=-\Delta $, the calculations give that the symbol-kernel of
$G^\pm(\log P)$ is \linebreak$\frac{-1}{x_n+y_n}e^{-|\xi
'|(x_n+y_n)}$ for $|\xi '|\ge 1$; the same holds for $P=\OP([\xi ]^2)$.
\endexample

When the order $m$ is even, there is a remarkable simplification in
view of Lemma 2.1:

\proclaim{Proposition 2.9} When $m=2k$, $k$ integer $>0$, then in
local coordinates, the symbol-kernel of $G^\pm(\log P)$ satisfies 
for $|\xi '|\ge 1:$
$$
\tilde g^{\pm}(\log p)(x',x_n,y_n,\xi ')=
\frac{-k}{x_n+y_n}e^{-|\xi
'|(x_n+y_n)}+\tilde g^{\pm, 0}(x',x_n,y_n,\xi '),\tag2.42
$$
where $\tilde g^{\pm, 0}(x',x_n,y_n,\xi ')$ is a standard singular
Green symbol of order and class $0$.

\endproclaim
\demo{Proof} We here have in view of Lemma 2.1 that the symbol of
$\log P$ is the sum of $k\log [\xi ]^2$ and a symbol $l(x,\xi )$ of
order 0 satisfying the transmission condition. Then we can apply Example 2.8 to
the first term and the standard $G^{\pm}$ construction (of \cite{G1})
to the second term.\qed
\enddemo

Thus in the even-order case, the terms in $G^{\pm}(\log P)$ of order $<0$ 
satisfy {\it all} the
standard s.g.o.\ estimates.

\subhead 3. Trace formulas
\endsubhead

The normal trace $\tr_n G$ of a singular Green operator $G$ with
symbol-kernel  $\tilde g(x',x_n,y_n,\xi ')$ in a local coordinate
system is the $\psi $do $S=\tr_n G$
with symbol 
$$
s(x',\xi ')=(\tr_n\tilde g)(x',\xi ')=\int_0^\infty \tilde g(x',x_n,x_n,\xi ')\, dx_n.\tag3.1
$$

In the differential operator case, we see from the estimates (2.23),
(2.28), (2.29) that $\tr_n\tilde g^{\log}_{-j}$ is well-defined for
$j\ge 1$. (Example 2.3 shows that this will generally not hold for the
principal 
part.) In view of the homogeneity (2.22), $\tr_n\tilde g^{\log}_{-j}$ is homogeneous of degree $-j$ in $\xi '$ for $|\xi '|\ge
1$, hence a classical $\psi $do symbol of degree $-j$. In the
pseudodifferential case, we have when $\nu >1$ and $j\ge 1$ that the
$L_{2,x_n,y_n}$-estimates  of $\tilde g^{\log}_{-j}$, $y_n\tilde
g^{\log}_{-j}$, $\partial_{y_n}\tilde g^{\log}_{-j}$ and
$y_n\partial_{y_n}\tilde g^{\log}_{-j}$ imply as in
in \cite{G1, pf.\ of Th.\ 3.3.9} that there is a well-defined normal trace, again a homogeneous classical
symbol of order $-j$. This estimation applies also to remainders  
$\tilde g^{\log}-\sum_{j<J}\tilde g^{\log}_{-j}$ for $J\ge 1$.

For $\nu =\frac12$ or $1$, the estimates in Theorem 2.6 do not provide
the estimates of $\partial_{y_n}\tilde g^{\log}_{-j}$ needed for this
argument. However, it is still posssible to take the normal trace of
$G_\lambda $, subtract the principal part, and integrate the remaining
operator with
$\log\lambda  $ to get a classical $\psi $do of order $-1$.

\proclaim{Theorem 3.1} In a local coordinate system, let $S_\lambda
=\tr_n G_\lambda$ with symbol $s(x',\xi ',\lambda )=(\tr_n\tilde
g)(x',\xi ',\lambda )$, expanded in terms $s_{-m-j}(x',\xi ',\lambda
)=(\tr_n\tilde g_{-m-j})(x',\xi ',\lambda )$. Define the parts
of $G_\lambda $ and $S_\lambda $ 
of order $-m-1$  by
$$
\aligned
G_ {\lambda ,\operatorname{sub}}&=G_\lambda - \operatorname{OPG}(\tilde
g_{-m}(x',x_n,y_n,\xi ',\lambda )) , \\
S_ {\lambda ,\operatorname{sub}}&=\tr_n G_ {\lambda ,\operatorname{sub}}=S_\lambda - \operatorname{OP}'(s_{-m}(x',\xi ',\lambda ))\\
\endaligned\tag3.2 
$$
(the remainders after subtracting principal parts),
and let 
$$ G^{\log}_{\operatorname{sub}}=\tfrac i{2\pi }\int_{\Cal C}\log \lambda
\,G_{\lambda ,\operatorname{sub}}\,d\lambda
,\tag3.3
$$
with symbol-kernel 
$\tilde g^{\log}_{\operatorname{sub}}=\tilde g^{\log}
-\tilde
g^{\log}_0$. The formula
$$
S^{\log }_{\operatorname{sub}}=\tfrac i{2\pi }\int_{\Cal C}\log \lambda
\,S_{\lambda,\operatorname{sub}} \,d\lambda 
\tag3.4
$$
defines a classical $\psi $do of order $-1$, with symbol
$s^{\log}_{\operatorname{sub}}(x',\xi ')$ expanded in terms
$$
s^{\log }_{\operatorname{sub},-j}(x',\xi ')=\tfrac i{2\pi }\int_{\Cal
C}\log \lambda \,s_{-m-j}(x',\xi ',\lambda ) \,d\lambda ,\quad j\ge 1.
\tag3.5
$$

When $\nu >1$, $S^{\log }_{\operatorname{sub}}$
 is the normal trace of $G^{\log}_{\operatorname{sub}}$.

\endproclaim

\demo{Proof} Since $G_\lambda $ and $G_{\lambda ,\operatorname{sub}}$ are of
regularity $\nu \ge \frac12$, $S_\lambda $ and $S_{\lambda ,\operatorname{sub}}$ are of
regularity $\nu -\frac14\ge \frac14$, cf.\ \cite{G3, Section 3}. In
particular, the symbols in $S_{\lambda ,\operatorname{sub}}$ are $O(\lambda ^{-1-1/4m})$
on the rays in $V$
so that the integrals in (3.4) and (3.5) make sense.

As accounted for in the text before the theorem, there are estimates
in the cases $\nu >1$ that allow interchange of the $\lambda
$-integral with the $x_n$-integral involved in taking $\tr_n$.\qed   
\enddemo

For the operator in Example 2.3, we note that 
$S_ \lambda =\tr _nG_\lambda $ is
the $\psi $do with symbol $-(2\kappa _1)^{-2}=-\frac14(\ang{ \xi
'}^2-\lambda )^{-1}$, so its log-integral gives $-\frac14\log(1- \Delta
_{x'})$. This demonstrates that the ``log-transform'' of the principal
part of $S_\lambda $ will not in general be a classical $\psi $do.

\medskip
Finally, we shall connect this with the study of the expansion
coefficient $C_0(I,(P+G)_T)$ 
in the last section of \cite{G3}; we here
write it simply as $C_0((P+G)_T)$ (or $C_0(B)$).
It is known from \cite{G1, Sect.\ 3.3} that when $m>n$, the trace of
the resolvent has an expansion in powers of $-\lambda $, 
$$
\Tr R_\lambda =\sum_{0\le l\le
n}c_{l}(-\lambda )^{\frac{n-l}m-1}+O(\lambda ^{-1-\frac1{4m}}),\tag3.6
$$
and a similar
proof shows that for general $m>0$, the expansion holds for a
sufficiently high iterate: 
$$
\Tr R_\lambda^N= \Tr \deln R_\lambda
 =\sum_{0\le l\le
n}c^{(N)}_{l}(-\lambda )^{\frac{n-l}m-N}+O(\lambda ^{-N-\frac1{4m}}).\tag3.7
$$
Define {\it the basic zeta value} as the coefficient of $(-\lambda )^{-N}$:
$$
C_0(B)=c ^{(N)}_n,
\tag3.8
$$
it is independent of $N$. If $B$ is invertible, $C_0(B)$
equals the value of the zeta function $ \zeta (B,s)$ --- the
meromorphic extension of $\Tr(B^{-s})$ --- at $s=0$. If $B$
has a nontrivial nullspace, the constants are connected by
$$
C_0(B)=\zeta (B,0)+\nu _0,
\tag3.9
$$
where $\nu _0$ is the dimension of the generalized eigenspace of the
zero eigenvalue.

There are similar expansions as in (3.7) of the traces of the $\psi$do
iterates $Q_\lambda ^N$ on $\widetilde X$, truncated to $X$, that
follow from integration over $X$ of the diagonal kernel expansions, as
established in \cite{G1, Sect.\ 3.3} (with remarks); it is the s.g.o.\
contribution that presents the greater challenge in \cite{G1}. In view
of the identifications in \cite{G3, Sect.\ 1}, the coefficient
of $(-\lambda )^{-N}$ here equals 
$-\tfrac1m\operatorname{res}_+(\log P)$, where the plus-index
indicates that  the pointwise
contribution to $-\tfrac1m\operatorname{res}(\log P)$ is integrated
over $X$ only. It can also be regarded as
$-\tfrac1m\operatorname{res}((\log  P)_+)$, extending the notation 
of \cite{FGLS}. 

The constant $C_0(B)$ was analyzed in \cite{G3, Sect.\ 5} in
relation to residue formulas, and we can now improve the result with
further information.

\proclaim{Theorem 3.2} One has that
$$
C_0(B)=-\tfrac1m\operatorname{res}_+(\log P)
-\tfrac1m\operatorname{res}_{X'}(S^{\log}_{\operatorname{sub}}),
\tag3.10
$$
where the terms are calculated as sums of contributions from local
coordinate patches of the form
$$
\int_{{\Bbb R}^n_+}\int_{|\xi |=1}\tr l_{-n}(x,\xi )\,\d
S(\xi)dx,\text{ resp. } \int_{{\Bbb R}^{n-1}}\int_{|\xi '|=1}\tr
s^{\log}_{\operatorname{sub},1-n} (x',\xi ')\,\d S(\xi ')dx'.
\tag3.11
$$
The term $-\tfrac1m\operatorname{res}_+(\log P)$ has an invariant
meaning as the coefficient of $(-\lambda )^{-N}$ in the expansion
similar to {\rm (3.7)} of 
$\Tr(((P-\lambda )^{-N})_+)$, 
and hence the last term in the right-hand
side of {\rm (3.10)} likewise has an invariant meaning.

When the problem is differential, or when the problem is
pseudodifferential with regularity $\nu >1$, then
$\operatorname{res}_{X'}(S^{\log}_{\operatorname{sub}})$ is, in local
coordinates, the residue of the normal trace of
$G^{\log}_{\operatorname{sub}}$.
\endproclaim

\demo{Proof}
It was shown in \cite{G3, Sect.\ 5} how $C_0(B)$ is found from
integrals of the strictly homogeneous symbol terms of order $-m-n$ in
$(P- \lambda )^{-1}$ resp.\ of order $-m-n+1$ in $G_\lambda $;
the proof given for the case $m>n$ extends to general $m$ when the
iterates are used, cf.\ \cite{G3, Remark 3.12}. It was shown moreover
that these integrals by use of \cite{G3, Lemmas 1.2, 1.3} could be
turned into log-integrals as   in (3.5). In those proofs, the
log-integration is applied after the $\tr_n$-integration, so the
boundary term is really $\operatorname{res} (S^{\log}_{\operatorname{sub}})$, as
defined in Theorem 3.1. 

When $\nu >1$, in particular when the problem is
differential so that $\nu =\infty $, Theorem {\rm 3.1} shows that
$S^{\log}_{\operatorname{sub}}$ is the normal trace of
$G^{\log}_{\operatorname{sub}}$, so the 
assertion for the residues follows.
\qed 
\enddemo

What we gain here in comparison with \cite{G3, Sect.\ 5} is a little more
insight into how the boundary term stems from the s.g.o.-like part of
$\log B$, plus the inclusion of all orders $m>0$. 
At any rate, since $C_0(B)$ 
is an invariant, we can propose it to
be the residue of $-\frac1m \log B$:

\proclaim{Definition 3.3} When $\{P_++G-\lambda ,T\}$ satisfies the
hypotheses of parameter-ellipticity given above, the residue of 
$\log(P+G)_T$ is defined to be the constant $$
\res(\log (P+G)_T)=-mC_0((P+G)_T)=\operatorname{res}_+(\log P)
+\operatorname{res}_{X'}(S^{\log}_{\operatorname{sub}}),\tag3.12
$$
as calculated in Theorem {\rm 3.2}. 
\endproclaim

This is consistent with the definition of \cite{FGLS}. We note
that certain steps in an explicit calculation of 
this constant depend very much on
localizations, e.g.\ in the steps of discarding the principal symbol and 
taking $\tr_n$. 
A number of similar or more general residue definitions are made in
\cite{G4} for compositions of $\psi $dbo's with components of $\log
P_T$ (when $P_T$ is defined from an
even-order differential problem). 
These residues do have a certain amount of
traciality: $
\res ([A,\log P_T])=0$ holds for operators $A$ of
order and class zero (cf.\ Theorem 6.5 there).

It should be noted that Definition 3.3 does not cover the case of
first-order differential 
operators with spectral boundary conditions, since such boundary
conditions are not {\it normal}. But for such boundary problems
(Atiyah-Patodi-Singer problems \cite{APS}) there
exists a wealth of other treatments, adapted to the specific
situation. The results there often depend on additional symmetry
properties. (See e.g.\ \cite{G2} and its references.)

\subhead 4. Sectorial projections
\endsubhead

Now we turn our attention to a certain spectral projection connected to
the realization $(P+G)_T$; namely a projection whose range contains the closure of the
direct sum of the generalized eigenspaces for the eigenvalues in a
sector of the complex plane. Such projections have been studied
earlier by Burak \cite{Bu}, Wodzicki \cite{W2}, and Ponge \cite{P}; the
latter gives a detailed deduction of the basic properties in the case
of classical $\psi $do's on closed manifolds. We recall the properties
below, supplying them with some additional information.

In order to apply the techniques to different types of operators, we
first consider an abstract situation where $A$ denotes an unbounded,
densely defined, closed operator in a Hilbert space $H$. 
It is assumed to have the following properties:

$A$ has a resolvent set containing two sectors $V_\theta$ and
$V_\varphi$ around $e^{i\theta} {\Bbb R}_+$ and $e^{i\varphi} 
{\Bbb R}_+$, respectively, for some $\theta < \varphi < \theta +
2\pi$, the resolvent $(A-\lambda )^{-1}$ is compact, and $\|
(A-\lambda)^{-1} \|$ is $O(\lambda^{-1})$ for $\lambda$ going to
infinity on each ray of these sectors. (We refer to Kato \cite{K} for
general background theory.)

For $x\in D(A)$ and $\lambda$ on a ray in either sector, we have
$$
\| \lambda^{-1} A\, (A-\lambda)^{-1} \, x \| \le \| \lambda^{-1}
(A-\lambda)^{-1} \| \cdot \| A\, x \| = O(\lambda^{-2}),
\tag 4.1
$$
so that $\lambda^{-1} A\, (A-\lambda)^{-1} \, x$ is integrable for
$|\lambda|\to \infty$.

Then define the operator $\Pi_{\theta,\varphi}(A)$, the {\it sectorial
projection}, with domain $D(A)$ to begin with, by 
$$
\Pi_{\theta,\varphi}(A) x = \tfrac{i}{2\pi} \int_{\Gamma_{\theta,\varphi}}
\lambda^{-1} A \, (A-\lambda)^{-1} \, x \, d\lambda, \quad x\in D(A),
\tag 4.2
$$
where the integration goes along the sectorial contour 
$$
\Gamma_{\theta,\varphi} = \{ r e^{i \varphi} \mid \infty > r > r_0
\} \cup \{ r_0 e^{i\omega} \mid \varphi \ge \omega \ge \theta
\} \cup \{ r e^{i \theta} \mid r_0 < r < \infty \},
\tag 4.3
$$
with $r_0$ taken so small that 0 is the only possible eigenvalue in
$\{|\lambda |<r_0\}$. If the operator is bounded in $H$-norm, we
extend it to $H$. 
This operator is a spectral projection in the following sense: 

For each $\lambda \in \sigma(A)$, denote the generalized eigenspace by
$E_\lambda $, 
$$
E_\lambda = \bigcup_{k \in \Bbb N} \ker (A-\lambda)^{k}
$$
(it equals $\ker (A-\lambda)^{k_0}$ for a sufficiently large $k_0$).
For $\alpha < \beta$, set 
$$
\Lambda_{\alpha,\beta} = \{\, r e^{i\omega}
\mid r>0,\,\alpha < \omega < \beta, \,\}, \quad E_{\alpha,\beta} =
\dotplus_{\lambda \in \sigma(A) \cap \Lambda_{\alpha,\beta}}
E_\lambda.
$$

\proclaim{Proposition 4.1} $\Pi_{\theta,\varphi}(A)^2 =
\Pi_{\theta,\varphi}(A)$, i.e.\ $\Pi_{\theta,\varphi}(A)$ is a
(possibly unbounded) projection in $H$. Its range contains
$E_{\theta,\varphi}$ and its kernel contains $E_0 \dotplus
E_{\varphi,\theta+2\pi}$.

{\rm (a)} If $A$ has a complete system of root vectors, i.e.\
$\dotplus_{\lambda\in \sigma(A)} E_\lambda$ is dense in $H$, then
$\Pi_{\theta,\varphi}(A)$ is the bounded projection onto
$\overline{E_{\theta,\varphi}}$  along $E_0 \dotplus \overline{
E_{\varphi,\theta+2\pi}}$.

{\rm (b)} If $A$ is normal, i.e.\ $A^* A = A A^*$, then
$\Pi_{\theta,\varphi}(A)$ is the bounded 
orthogonal projection onto $\oplus_{\lambda\in\sigma(A) \cap
\Lambda_{\theta,\varphi}} \ker(A-\lambda)$ along
$\oplus_{\lambda\in\sigma(A) \setminus \Lambda_{\theta,\varphi}}
\ker(A-\lambda)$.
\endproclaim

\demo{Proof}
Except for a few elementary considerations regarding the domain and
closedness, the proofs of \cite{P, Propositions 3.2, A.4, and A.5} carry
over almost word for word to the present setting (it should be noted
that some contours in \cite{P} have the opposite orientation).

In (a) and (b), the boundedness of $\Pi_{\theta,\varphi}(A)$ follows
from the fact that the kernel and range are closed.
\qed
\enddemo

In certain important cases, $\Pi_{\theta,\varphi}(A)$ can be seen to
be bounded regardless of whether the hypotheses of (a) or (b) can be
verified; as shown in \cite{P, Proposition 3.1} this holds when $A$
is a $\psi$do of order $m>0$ on a closed manifold. We shall see below in
Theorem 4.6 that it also holds for the realization of a differential
elliptic boundary value problem.

As shown below, the sectorial projection has a direct connection
with the choice of spectral cut in our definition of the logarithm of
an operator. Using arguments as in Section 2, we can define the
logarithm of $A$ with a branch cut at the angle $\theta$ as
$$
\log_\theta \! A = \lim_{s\searrow 0} \tfrac{i}{2\pi}
\int_{\Cal C_\theta} \lambda_\theta^{-s} \log_\theta\! \lambda \,
(A-\lambda)^{-1} \, d\lambda
\tag 4.4
$$
where the subscript $\theta$ indicates that $\lambda^{-s} \log
\lambda$ is chosen to have a branch cut along
$e^{i\theta}\Bbb R_+$, and the contour is the Laurent loop
$$
\Cal C_\theta=\{ re^{i\theta }\mid \infty >r>r_0\}\cup\{ r_0e^{i\omega
}\mid \theta \ge \omega \ge \theta - 2\pi 
\}\cup\{re^{i (\theta-2\pi) }\mid r_0<r<\infty \}.
\tag 4.5
$$

The following proposition eliminates the limiting procedure of (4.4)
and gives a useful alternative description of
$\Pi_{\theta,\varphi}(A)$. A proof can be found in the Appendix.

\proclaim{Proposition 4.2}
For $x \in D(A)$ we have the identities
$$
\align
\log_{\theta} \! A \, x & = \tfrac{i}{2\pi} \int_{\Cal C_{\theta}}
\lambda^{-1} \log_\theta \! \lambda \, A (A-\lambda)^{-1} \, x \,
d\lambda \quad \text{ and }
\tag 4.6 \\
\Pi_{\theta,\varphi}(A) \, x & = \tfrac{i}{2\pi}
\int_{\Gamma_{\theta,\varphi}} (A-\lambda)^{-1} \, x \, d\lambda +
\frac{\varphi-\theta}{2\pi} \, x,
\tag {4.7}
\endalign
$$
where the integral in the right-hand side of {\rm (4.7)} is an
improper integral. 
\endproclaim

Next, we include a lemma which will be useful for our considerations
regarding expressions involving different branches of the
logarithm. Again, a proof is available in the Appendix.

\proclaim{Lemma 4.3}
Let $f(\lambda)$ be a continuous (possibly vector-valued) function on
the ``punctuated double keyhole region'' 
$$
V_{r_0,\delta} = \{ \lambda \in \Bbb C \mid |\lambda| < 2r_0 \text{ or
} |\arg \lambda - \theta| < \delta \text{ or } |\arg\lambda -
\varphi|< \delta \} \setminus \{ 0 \},
\tag 4.8
$$
such that $f(\lambda)$ is $O(\lambda^{-1-\varepsilon})$ for $|\lambda|
\to \infty$ in $V_{r_0,\delta}$. Then
$$
\int_{\Cal C_{\theta}} \log_\theta \! \lambda \,
f(\lambda) \,d\lambda - \int_{\Cal C_{\varphi}} \log_\varphi \! \lambda \,
f(\lambda)\, d\lambda = -2\pi i \int_{\Gamma_{\theta,\varphi}} f(\lambda)
\,d\lambda.
\tag 4.9
$$
\endproclaim

We can use this lemma to describe the relation between
$\Pi_{\theta,\varphi}(A)$ and logarithms of $A$ as follows:

\proclaim{Proposition 4.4}
For $x \in D(A)$,
$$
\log_\theta \! A \, x - \log_\varphi \! A \, x =
\int_{\Gamma_{\theta,\varphi}} \lambda^{-1} A (A-\lambda)^{-1} \, x \,
d\lambda = -2\pi i \; \Pi_{\theta,\varphi}(A) \, x.
\tag 4.10
$$
When $\Pi_{\theta,\varphi}(A)$ is bounded, so is
$\log_\theta \! A  - \log_\varphi \! A  $, and 
$$
\Pi_{\theta,\varphi}(A)=\tfrac i{2\pi }( \log_\theta \! A -
\log_\varphi \! A ).
\tag 4.11
$$
\endproclaim

\demo{Proof} For $x\in D(A)$, the expression $f(\lambda) = \lambda^{-1}
A (A - \lambda)^{-1} x$ is holomorphic in
$V_{r_0,\delta}$ for some $r_0, \delta > 0$, and $f(\lambda)$ is
$O(\lambda^{-2})$ for $|\lambda|\to \infty$ in $V_{r_0,\delta}$ by (4.1).

Hence we can apply Lemma 4.3, and insertion of the expression for
$f(\lambda)$ into (4.9) gives
$$
\multline
\int_{\Cal C_{\theta}} \!\! \log_\theta \! \lambda \,
\lambda^{-1} A (A-\lambda)^{-1} \, x \, d\lambda - \int_{\Cal
C_{\varphi}} \!\! \log_\varphi \! \lambda \, 
\lambda^{-1} A (A-\lambda)^{-1} \, x \, d\lambda \\
= -2\pi i \int_{\Gamma_{\theta,\varphi}} \!\! \lambda^{-1} A
(A-\lambda)^{-1} \, x \, d\lambda.
\endmultline
\tag 4.12
$$
Then (4.10) follows from (4.2) and (4.6).

If $\Pi_{\theta,\varphi}(A)$ is bounded, (4.10) extends to all
$x\in H$ since $D(A)$ is dense in $H$, and (4.11) follows.
\qed
\enddemo

With the results above at hand we return to the realization
$(P+G)_T$. Modifying the assumption of Section 2 a little, we now
assume $\{P_+ + G - \lambda, T\}$ to satisfy the conditions of
parameter-ellipticity in \cite{G1, Def. 3.3.1} for $\lambda$ on the
rays of {\it two} sectors around $e^{i\theta} \Bbb R_+$ and
$e^{i\varphi} \Bbb 
R_+$, respectively. Then the realization $B=(P+G)_T$ satisfies the
requirements for $A$ given above, and we can define the sectorial
projection accordingly:
$$
\Pi_{\theta,\varphi}(B) = \tfrac{i}{2\pi}
\int_{\Gamma_{\theta,\varphi}} \lambda^{-1} B \, R_\lambda \,
d\lambda.
\tag 4.13
$$
Like in the case of the logarithm, we decompose it into the
contributions from the pseudo\-differential and singular Green parts.

For the $\psi$do $P$ on the closed
manifold $\widetilde X$, we can use Proposition 4.2 to see that
$$
 \frac{i}{2\pi} \int_{\Gamma_{\theta,\varphi}} Q_{\lambda}u \,
d\lambda + \frac{\varphi-\theta}{2\pi}u =
\Pi_{\theta,\varphi}(P)u,\quad u\in D(P);
\tag 4.14
$$
it is known from \cite{W2}, \cite{P}, that  $\Pi_{\theta,\varphi}(P)$ is
a $\psi$do of order
$\le 0$ on $\widetilde X$. 

Using Proposition 4.2, (2.4), and the fact that $r^+ e^+ = I$, we can
rewrite (4.13) as 
$$
\aligned
\Pi_{\theta,\varphi}(B) & = \frac{i}{2\pi}
\int_{\Gamma_{\theta,\varphi}} \, R_\lambda \, d\lambda +
\frac{\varphi-\theta}{2\pi} = \frac{i}{2\pi}
\int_{\Gamma_{\theta,\varphi}} [ Q_{\lambda,+} + G_\lambda ] \, d\lambda +
\frac{\varphi-\theta}{2\pi} \\
& = r^+ \big( \frac{i}{2\pi} \int_{\Gamma_{\theta,\varphi}} Q_{\lambda} \,
d\lambda + \frac{\varphi-\theta}{2\pi} \big) e^+ + \frac{i}{2\pi}
\int_{\Gamma_{\theta,\varphi}} G_\lambda \, d\lambda\\
& = \Pi_{\theta,\varphi}(P)_+ +  \frac{i}{2\pi}
\int_{\Gamma_{\theta,\varphi}} G_\lambda \, d\lambda;
\endaligned
\tag 4.15
$$
in the last line we moreover used (4.14).
Now an application of Proposition 4.4 to $P$ and $B$ gives:
$$
\aligned
\Pi_{\theta,\varphi}(P)_+ & = \tfrac{i}{2\pi} \big( (\log_\theta \!
P)_+ - (\log_\varphi \! P)_+ \big), \\
\Pi_{\theta,\varphi}(B) & = \tfrac{i}{2\pi} \big( \log_\theta \! B -
\log_\varphi \! B \big).
\endaligned
\tag 4.16
$$

Using the contour $\Cal C_\theta$ from (4.5) we can define an
operator as in (2.13),
$$
G^{\log_\theta} = \tfrac{i}{2\pi} \int_{\Cal C_\theta} \log _\theta \!\lambda \,
G_\lambda \, d\lambda,
\tag 4.17
$$
and similarly define $G^{\log_\varphi}$ where $\theta$ is replaced by
$\varphi$. By rotation it is obvious that $G^{\log_\theta}$ and
$G^{\log_\varphi}$ have properties similar to those of $G^{\log}$
described in Section 2. Now (4.16) and (2.15) show that if we define
$G_{\theta,\varphi }$ by
$$
G_{\theta,\varphi} = \tfrac{i}{2\pi} \int_{\Gamma_{\theta,\varphi}}
G_\lambda \, d\lambda ,\tag 4.18
$$
then
$$
G_{\theta,\varphi} = \tfrac{i}{2\pi} \big( G^{\log_\theta} -
G^{\log_\varphi} \big).
\tag 4.19
$$

In view of (4.15), we have then obtained:

\proclaim{Theorem 4.5} The sectorial projection for $B=(P+G)_T$ satisfies
$$
\Pi_{\theta,\varphi}(B) = \Pi_{\theta,\varphi}(P)_+ + G_{\theta,\varphi},
\tag 4.20
$$
where each term on the right hand side is known: $\Pi_{\theta,\varphi}(P)_+$
is the truncation of a $\psi$do on $\widetilde X$ of order at most
zero, in particular it is bounded on $L_2(X,E)$; $G_{\theta,\varphi}$ is
a difference {\rm (4.19)} of two terms  of the log-type described in
Section {\rm 2} and hence is a generalized singular Green operator,
bounded from $L_2(X,E)$ to $H^{-\varepsilon}(X,E)$.
\endproclaim

Like $G^{\log}$, $G_{\theta,\varphi}$ acts as in (2.21). It has a
symbol-kernel $\tilde g_{\theta,\varphi} \sim \sum_{j\in\Bbb N} \tilde
g_{\theta,\varphi,-j}$, with terms given by
$$
\tilde g_{\theta,\varphi,-j} = \tfrac{i}{2\pi} ( \tilde
g^{\log_\theta}_{-j} - \tilde g^{\log_\varphi}_{-j} ) = \tfrac{-1}{4\pi^2}
\big(
\int_{\Cal C_{\theta}} \log_\theta\! \lambda \, \tilde g_{-m-j}\, d\lambda -
\int_{\Cal C_{\varphi}} \log_\varphi\! \lambda \, \tilde g_{-m-j}\, d\lambda
\big).
\tag 4.21
$$
By Lemma 4.3 this is simplified to
$$
\tilde g_{\theta,\varphi,-j}(x',x_n,y_n,\xi') = \tfrac{i}{2\pi}
\int_{\Gamma_{\theta,\varphi}} \tilde g_{-m-j}(x',x_n,y_n,\xi',\lambda)
\,d\lambda.
\tag 4.22
$$

In view of (4.19) and (4.21), the results on $G^{\log}$ resp.\ $\tilde
g^{\log}$ in  Section 2 carry over immediately to $G_{\theta,\varphi}$ resp.\
$\tilde g_{\theta,\varphi}$. We shall not reproduce
all the statements explicitly, but will just present 
the following important result obtained from Theorem 2.4.

\proclaim{Theorem 4.6} Assume that $P$ is a differential operator,
$G=0$, and the trace operators $T_0,\ldots,T_{m-1}$ are differential
operators; hereby $B=P_T$.

Then $G_{\theta,\varphi}$ is, in local coordinates near $X'$, a
generalized singular Green operator 
$$
G_{\theta,\varphi} = \operatorname{OPG}(\tilde g_{\theta,\varphi})
\tag 4.23
$$ 
with $\tilde g_{\theta,\varphi} \sim \sum_{j\in\Bbb N} \tilde
g_{\theta,\varphi,-j}$; the $j$'th term is quasihomogeneous as in
{\rm (2.22}) and satisfies estimates as in {\rm (2.23)}.

$G_{\theta,\varphi}$ and $\Pi_{\theta,\varphi}(P_T)$ are
bounded operators in $L_p(X,E)$ for $1<p<\infty$. In particular,
$\Pi_{\theta,\varphi}(P_T)$ is a bounded projection in $L_2(X,E)$.
\endproclaim

\demo{Proof}
The claims regarding $\tilde g_{\theta,\varphi}$ follow immediately from
Theorem 2.4 and (4.21).

The boundedness properties of $G_{\theta,\varphi}$ are obvious from
Theorem 2.4 and (4.19). Since $\Pi_{\theta,\varphi}(P)_+$ is the
truncation of a $\psi$do of order at most zero, this is also bounded
in $L_p(X,E)$; then in view of (4.20) so is $\Pi_{\theta,\varphi}(P_T)$.
\qed
\enddemo

An interesting question is whether one can give criteria on $P$, $G$,
and $T$ assuring that the operator $\Pi_{\theta,\varphi}((P+G)_T)$
belongs to the Boutet de Monvel calculus.

Concerning the $\psi $do part $\Pi _{\theta ,\varphi }(P)$, with
symbol $\pi _{\theta ,\varphi }(x,\xi )$ in local coordinates, we have
easily by use of Lemma 2.1:

\proclaim{Lemma 4.7} When $m$ is even, $\pi _{\theta ,\varphi }(x,\xi
)$ satisfies the transmission condition.

Hence $\Pi_{\theta,\varphi}(P)_+$ is in the Boutet de Monvel calculus
for even $m$.
\endproclaim

\demo{Proof}  We have that in view of (2.10) that
$$
\operatorname{symb}(\log _\theta P)=m\log[\xi ]+l_\theta (x,\xi ),\quad
l_\theta (x,\xi )\sim\sum_{j\in {\Bbb N}}l_{\theta ,-j}(x,\xi ),
\tag 4.24
$$
where $m\log[\xi ]+l_{\theta ,0}(x,\xi )=\log _\theta (p_m(x,\xi ))$,
with similar formulas for $\log _\varphi P$, so 
the symbols of  $\log _\theta P$ and $\log _\varphi
P$ have the same log-term $m\log [\xi ]$.
Then it is seen from the first line in (4.16) that
$$
\pi _{\theta ,\varphi }(x,\xi )=\tfrac i{2\pi }(l_\theta (x,\xi
)-l_\varphi (x,\xi )),
\tag 4.25
$$
which satisfies the transmission condition when $m$ is even in view
of Lemma 2.1. \qed
\enddemo

This could also be based more directly on the fact, worked out in
detail in \cite{P}, that 
$\pi_{\theta ,\varphi}(x,\xi) \sim \sum_{j\in\Bbb N}
\pi_{\theta,\varphi,-j}(x,\xi)$, where the terms are given by
$$
\pi_{\theta,\varphi,-j}(x,\xi) = \tfrac{i}{2\pi} \int_{\Cal
C_{\theta,\varphi}(x,\xi)} q_{-m-j}(x,\xi,\lambda) \, d\lambda;
\tag 4.26
$$
here $\Cal C_{\theta,\varphi}(x,\xi)$ is a closed curve in the sector
$\Lambda_{\theta,\varphi}$ going in the positive direction around the
part of the spectrum of $p_m(x,\xi)$ lying in that sector.

When $m$ is odd, one cannot expect $\Pi_{\theta,\varphi}(P)$ to
satisfy the transmission condition. For example, for a first-order
selfadjoint invertible elliptic differential operator $A$ on
$\widetilde X$ (e.g., a Dirac operator), $\Pi _{-\frac{\pi
}2,\frac{\pi }2}(A)$ equals $\Pi _>(A)$, the positive
eigenprojection $\frac12(I+A|A^2|^{-1/2})$, where $A|A^2|^{-1/2}$ {\it
does not} satisfy the transmission condition (its even-order symbol
terms are odd in $\xi $).

Next, let us consider the s.g.o.\ part $G_{\theta,\varphi}$. Example
4.8 below shows a differential operator realization 
where $G_{\theta ,\varphi }$ is not a standard
singular Green operator, already in a constant-coefficient principal
symbol case. Example 4.9 on the other hand defines a
general class of differential operator realizations where
$G_{\theta,\varphi}$ is a standard s.g.o.\, and $\Pi
_{\theta,\varphi}(B)$ belongs to the standard 
calculus. Here one finds however, that lower order perturbations can
ruin the standard s.g.o.-properties.

\example{Example 4.8} Consider the differential operators $A$ and
$P$ on $\Bbb R^4_+$ given by
$$
A = \left( \matrix i & 0 \\ 0 & -i \endmatrix \right) D_1 +
\left( \matrix 0 & 1 \\ -1 & 0 \endmatrix \right) D_2 + \left( \matrix 0
  & i \\ i & 0 \endmatrix \right) D_3 + \left( \matrix 1 & 0 \\ 0 & 1
\endmatrix \right) D_4,
\tag 4.27
$$
and 
$$
P = \left( \matrix 0 & -A^* \\ A & 0 \endmatrix \right),
\tag 4.28
$$
where $A^*$ denotes the formal adjoint of $A$. ($A$ and $P$ are
Dirac-type operators, with $A^*A=-\Delta  I_2$, $(iP)^2=-\Delta I_4$.)

Regarding this as a localization of a manifold situation, we seek the
projection onto the (generalized) eigenspaces for the eigenvalues
$\lambda$ in the upper halfplane $\Bbb C_+$ for a certain realisation
$P_T$ of $P$, where the boundary condition is $B\gamma_0u=0$, with
$$
B = \left( \matrix 1 & 0 & 1 & 0 \\ 0 & 1 & 0 & 1 \endmatrix \right),
\tag 4.29
$$
i.e., $\gamma_0 u_1 + \gamma_0 u_3 = \gamma_0 u_2 + \gamma_0 u_4 = 0$,
$u_i$ being the $i$'th component of $u$.

Thus, in this localized situation we shall construct
$\Pi_{\theta,\varphi}(P_T)$ with $\theta = 0$ and $\varphi = \pi$. In
this case the contour $\Gamma_{\theta,\varphi}$ is a contour from
$-\infty$ to $\infty$ passing above the origin.

$P$ has symbol
$$
p(\xi) = \left( \matrix 0 & -{}^t\overline{a(\xi)}
  \\ a(\xi) & 0 \endmatrix \right) =  \left( \matrix 0 & 0 & i\xi_1
  - \xi_4 & \xi_2 + i \xi_3 \\ 0 & 0 & -\xi_2 + i \xi_3 & -i\xi_1 -
  \xi_4 \\  i\xi_1 + \xi_4 & \xi_2 + i \xi_3 & 0 & 0 \\
  -\xi_2 + i \xi_3 & -i\xi_1 + \xi_4 & 0 & 0 \endmatrix \right),
$$
the eigenvalues of which are $\pm i |\xi|$. Hence $P-\lambda$ is
parameter-elliptic for $\lambda$ on all rays in $\Bbb C\setminus i
\Bbb R$, with parametrix-symbol
$$
q(\xi,\lambda) = (p(\xi) - \lambda)^{-1} = \frac{1}{|\xi|^2 +
  \lambda^2} \left( \matrix -\lambda & 0 & -i\xi_1
  +\xi_4 & -\xi_2 - i \xi_3 \\ 0 & -\lambda & \xi_2 - i \xi_3 & i\xi_1 +
  \xi_4 \\  -i\xi_1 - \xi_4 & -\xi_2 - i \xi_3 & -\lambda & 0 \\
  \xi_2 - i \xi_3 & i\xi_1 - \xi_4 & 0 & -\lambda \endmatrix \right).  
$$
We first find the $\psi$do part of $\Pi_{0,\pi}(P_T)$: According to
(4.26) the symbol $\pi(\xi)$ of
$\Pi_{0,\pi}(P)$ is obtained by integrating $q(\xi,\lambda)$ along a
small closed curve, $\Cal C_\xi$, enclosing the pole $i 
|\xi|$ in $\Bbb C_+$:
$$
\multline
\pi(\xi) = \tfrac{i}{2\pi} \int_{\Cal C_\xi} q(\xi,\lambda) d\lambda =
- \Res_{\lambda = i |\xi|} \big( q(\xi,\lambda) \big) \\ =  \frac{1}{2|\xi|} \left(
  \matrix |\xi| & 0 & \xi_1 + i\xi_4 & -i \xi_2 + \xi_3 \\
0 & |\xi| & i \xi_2 + \xi_3 & -\xi_1 + i \xi_4 \\
\xi_1 - i \xi_4 & -i\xi_2 + \xi_3 & |\xi| & 0 \\
i\xi_2 + \xi_3 & - \xi_1 - i\xi_4 & 0 & |\xi| \endmatrix \right).
\endmultline
\tag 4.30
$$

The singular Green part $G_\lambda$ of the
resolvent $R_\lambda = (P_T - \lambda)^{-1}$ has symbol-kernel
$$
\tilde g(x_n,y_n,\xi',\lambda) = \frac{1}{2\sigma} \left( \matrix -i
\xi_1 + i\sigma & - \xi_2-i\xi_3 & - \lambda & 0 \\ \xi_2 -i \xi_3 &
i\xi_1 + i\sigma & 0 & -\lambda \\
-\lambda & 0 & -i\xi_1 -i \sigma & -\xi_2 - i \xi_3 \\
0 & - \lambda & \xi_2 - i\xi_3 & i\xi_1 - i \sigma \endmatrix \right)
e^{-\sigma(x_n+y_n)},
$$
where $\sigma = \sqrt{|\xi'|^2+\lambda^2}$. Note that $\sigma$ is
holomorphic (and $\operatorname{Re} \sigma > 0$) for $\lambda \in \Bbb C \setminus
\pm i (|\xi'|,\infty)$; in particular $\{P - \lambda, B\gamma _0\}$ is
parameter-elliptic for $\lambda$ on any ray in $\Bbb C\setminus i\Bbb R$.

The integration contour $\Gamma_{0,\pi}$ is homotopic in $\{ r
e^{i\omega} \mid \omega \neq \pm \tfrac \pi2 \text{ or } r < |\xi'|
\}$ to the real line; thus, due to the exponential falloff of
$e^{-(|\xi'|^2+\lambda^2)^{\frac 12} (x_n+y_n)}$ we get
$$
\tilde g_{\theta,\varphi}(x_n,y_n,\xi') = \tfrac{i}{2\pi}
\int_{\Gamma_{\theta,\varphi}} \tilde g(x_n,y_n,\xi',\lambda) \, d\lambda
= \tfrac{i}{2\pi} \int_{-\infty}^\infty \tilde g(x_n,y_n,\xi',t)\, dt.
\tag 4.31
$$
We can now verify that $\tilde g_{\theta,\varphi}$ is not a singular
Green symbol-kernel: The 12-matrix entry of $\tilde
g_{\theta,\varphi}$ becomes
$$
\frac{-i\xi_2 +\xi_3}{4\pi} \int_{-\infty}^\infty
(|\xi'|^2+t^2)^{-\frac 12} e^{-(|\xi'|^2+t^2)^{\frac 12} (x_n+y_n)} \,dt,
\tag4.32$$
which, for fixed $\xi'$, is unbounded as $x_n + y_n$ goes to zero;
hence, $\tilde g_{\theta,\varphi}$ is not in $\Cal S_{++}$.

To see this note that, for fixed $a>0$,
$$
\aligned
f(r) & = \tfrac 12 \int_{-\infty}^\infty (a^2+t^2)^{-\frac 12} e^{-r
(a^2+t^2)^{\frac 12}} dt = \int_0^\infty (a^2+t^2)^{-\frac 12} e^{-r
(a^2+t^2)^{\frac 12}} dt \\
& \ge \int_0^\infty \frac{e^{-(a+t) r }}{a+t} dt = \int_{a r}^\infty
\frac{e^{-u}}{u} du 
\endaligned
$$
which diverges to $+\infty$ as $r\to 0^+$.

\endexample

\example{Example 4.9}
Let $X'_0$ be a closed $(n-1)$-dimensional manifold provided with an
elliptic second-order differential operator $S$ which is selfadjoint
positive in $L_2(X'_0)$. Let $X=X'_0\times [0,a]$ with 
points $x=(x',x_n)$,
$x'\in X'_0$ and $x_n\in [0,a]$,
and let $B$ be the Dirichlet
realization of  $D_{x_n}^2+S$ on $X$; it is  selfadjoint positive in 
$L_2(X)$, with $D(B)=H^2(X)\cap H^1_0(X)$. Let $A$ be the Dirichlet
realization of   
$$
P = \left( \matrix D_{x_n}^2+S & S \\ S & -D_{x_n}^2-S \endmatrix \right).
\tag 4.33
$$
on $X$, then in
fact, 
$$
A=\pmatrix B& S\\ S&-B\endpmatrix \tag4.34
$$
with domain $D(B)\times D(B)$. The resolvent is 
$$
(A-\lambda )^{-1}=\pmatrix -B-\lambda & -S\\ -S&B-\lambda \endpmatrix
(\lambda ^2-B^2-S^2)^{-1},\tag4.35
$$
where we used that $S$ and $B$ commute.
Define $B_1=(B^2+S^2)^\frac12$. Here $B^2+S^2$ is the realization
of the fourth-order
elliptic differential operator $(D_{x_n}^2+S)^2+S^2$ determined by 
the boundary condition
$\gamma _0u=0, \gamma _0Bu=0$. This is one of the particular cases where the
square root of the interior operator does satisfy the transmission
condition,
cf.\ \cite{G1, (4.4.9)}. Moreover, the square root of the realization
$B^2+S^2$ represents a boundary condition consisting of exactly the
part of the boundary condition for $B^2+S^2$ that makes sense on
$H^2(X)$, cf.\ \cite{G1, Cor.\ 4.4.3} (based on a result of Grisvard);
so in fact $B_1$ is the realization of $((D_{x_n}^2+S)^2+S^2)^\frac12$
determined by the Dirichlet condition $\gamma _0u=0$. This belongs to
the standard calculus and enters nicely in the theory of
\cite{G1}, cf.\ Section 1.7 there. Note that $D(B_1)=D(B)$.

We can then calculate
$$
\aligned
(\lambda ^2-(B^2+S^2))^{-1}&=(\lambda ^2-B_1^2)^{-1}
=(B_1-\lambda )^{-1}(-B_1-\lambda )^{-1}\\
&=(B_1-\lambda )^{-1}(2B_1)^{-1}(B_1+\lambda +B_1-\lambda )(-B_1-\lambda
)^{-1}\\
&=-\tfrac12 B_1 ^{-1}\bigl((B_1-\lambda )^{-1}-(-B_1-\lambda )^{-1}\bigr)
\endaligned,\tag4.36
$$
 which leads to the formula:
$$
\aligned
(A-\lambda )^{-1}&=\pmatrix -B+B_1-B_1-\lambda & -S\\
-S&B-B_1+B_1-\lambda  \endpmatrix(B_1-\lambda )^{-1}(-B_1-\lambda )^{-1}\\
&=
\pmatrix (B_1-\lambda )^{-1}&0\\0&(-B_1-\lambda )^{-1}\endpmatrix\\
&\quad-\pmatrix B_1-B & -S\\
-S&B-B_1  \endpmatrix\tfrac12 B_1 ^{-1}\bigl((B_1-\lambda
)^{-1}-(-B_1-\lambda )^{-1}\bigr), 
\endaligned\tag4.37
$$
valid for $\lambda $ outside the spectra of $B_1$ and $-B_1$.
%(We use that  $S$ and $B_1^{-1}$ commute with $(\pm B_1-\lambda
%)^{-1}$.) 
To determine the spectral projection $\Pi _{\theta ,\varphi }(A)$ with
$\theta =-\frac\pi 2$, $\varphi =\frac\pi 2$, we use the abstract
machinery. It is seen from either of the formulas (4.2) or (4.7) that
$$
\aligned
\Pi _{-\frac \pi 2,\frac\pi 2 }(A)&
=
\pmatrix \Pi _{-\frac \pi 2,\frac\pi 2 }(B_1)&0\\0&\Pi _{-\frac \pi 2,\frac\pi 2 }(-B_1 )\endpmatrix\\
&\quad-\pmatrix B_1-B & -S\\
-S&B-B_1  \endpmatrix\tfrac12 B_1 ^{-1}\bigl(\Pi _{-\frac \pi 2,\frac\pi 2
}(B_1)-\Pi _{-\frac \pi 2,\frac\pi 2 }(-B_1)\bigr). 
\endaligned\tag4.38$$
Here
$$
\Pi _{-\frac \pi 2,\frac\pi 2
}(B_1)=I ,\quad
\Pi _{-\frac \pi 2,\frac\pi 2
}(-B_1)=0 ,
\tag4.39
$$
in view of Proposition 4.1 and the fact that $B_1$ is selfadjoint
positive. It follows that$$
\Pi _{-\frac \pi 2,\frac\pi 2 }(A)
=\pmatrix \frac12+\frac12 B B_1 ^{-1}&\tfrac12 SB_1 ^{-1}\\\tfrac12
SB_1 ^{-1}& \frac12 -\frac12 BB_1 ^{-1}\endpmatrix.\tag4.40
$$
The operator is in the Boutet de Monvel calculus. Note that the sum of
the diagonal terms is $I$, so the residue of the operator is zero.

Inherent in this example are some symbol calculations where the poles
of the resolvent symbol appear isolated in such a way that integrals over
$\Gamma _{\theta ,\varphi }$ can be turned into integrals over
closed curves, reducing to simple residue
calculations. Perturbations can easily introduce more complicated
calculations where integrals as in (4.32) appear, leading to
non-standard s.g.o.-symbols (we shall not reproduce examples here).

\endexample

In view of Definition 3.3 and the formulas (4.16), the sectorial
projection $\Pi
_{\theta,\varphi}(B)$
has a well-defined residue. In the differential operator
case where the order $m$ is even, one can moreover define residues of the compositions of 
$\Pi
_{\theta,\varphi}(B)$ with operators $A$ in the Boutet de Monvel calculus;
this is taken up in \cite {G4}. It is found there that if in addition, $A$ is of order and class
0, the residue vanishes on the commutator of 
$\Pi
_{\theta,\varphi}(B)$
and $A$.

It is still an open question whether the residue is zero on sectorial
projections for boundary value problems, as it is in the closed
manifold case; we expect to return to this question in a
forthcoming work.

\subhead Appendix A. Proofs of auxiliary results in functional analysis 
\endsubhead

\demo{Proof of Proposition {\rm 4.2}}
First we prove (4.6): Let, for $N\in\Bbb N$, 
$$
  \Cal C_\theta^N = \{ r e^{i\theta} \mid N \ge r \ge r_0 \}
  \cup \{ r_0 e^{i \omega} \mid \theta \ge \omega \ge \theta-2\pi \}
  \cup \{ r e^{i (\theta - 2\pi)} \mid r_0 \le r \le N \}.
\tag A.1
$$
Then, for $s>0$,
$$
\int_{\Cal C_\theta^N} \lambda_\theta^{-s-1} \log_\theta\! \lambda
\,d\lambda = \big[ -\tfrac{1}{s^2} \lambda_\theta^{-s} ( 1 +
s\log_\theta \!\lambda ) \big]_{Ne^{i(\theta-2\pi)}}^{Ne^{i\theta}}
\longrightarrow 0 \text{ for } N\to\infty,
\tag A.2
$$
since $N^{-s}$ and $N^{-s} \log N$ go to $0$ for $N\to\infty$. It
follows that 
$$
\lim_{s\searrow 0} \lim_{N\to\infty} \int_{\Cal C_\theta^N}
\lambda_\theta^{-s-1} \log_\theta \!\lambda \,d\lambda = 0.
\tag A.3
$$
Observe that the order of the limits is important.

Using the resolvent identity $A(A-\lambda)^{-1} = 1 + \lambda
(A-\lambda)^{-1}$ we now get for $x\in D(A)$:
$$
\align
\lim_{s\searrow 0} \int_{\Cal C_\theta} \lambda_\theta^{-s}
&\log_\theta \! \lambda \, (A-\lambda)^{-1} \, x \, d\lambda 
= \lim_{s\searrow 0} \lim_{N\to\infty} \int_{\Cal C_\theta^N}
\lambda_\theta^{-s} \log_\theta \!\lambda \, (A-\lambda)^{-1} \, x \,
d\lambda \\
& = \lim_{s\searrow 0} \lim_{N\to\infty} \Big[ \int_{\Cal C_\theta^N}
\lambda_\theta^{-s-1} \log_\theta \!\lambda\, x \, d\lambda + \int_{\Cal
C_\theta^N} \lambda_\theta^{-s} \log_\theta \!\lambda \, (A-\lambda)^{-1}
\, x \, d\lambda \Big] \\
& = \lim_{s\searrow 0} \lim_{N\to\infty} \int_{\Cal C_\theta^N}
\lambda_\theta^{-s-1} \log_\theta\!\lambda\, \big[ 1 + \lambda
(A-\lambda)^{-1} \big] \, x \,
d\lambda \tag A.4 \\
& = \lim_{s\searrow 0} \lim_{N\to\infty} \int_{\Cal C_\theta^N}
\lambda_\theta^{-s-1} \log_\theta\!\lambda\, A (A-\lambda)^{-1} \, x \,
d\lambda,
\endalign
$$
where we used (A.3) in the second line (adding zero).
Then, since $\| (A-\lambda)^{-1} \| \, \dot\le \, |\lambda|^{-1}$,
$$
\aligned
  \| \lambda_\theta^{-s-1} \log_\theta \!\lambda \, A (A-\lambda)^{-1} x
\| \ \dot\le \ |\log \lambda\, | |\lambda|^{-s-2} \| A x\|,
\endaligned
\tag A.5
$$
so that the integrand in the last expression of (A.4) is integrable
along $\Cal C_\theta$ uniformly in $s>0$, and
$$
  \lim_{s\searrow 0} \lim_{N\to\infty} \int_{\Cal C_\theta^N}
    \lambda_\theta^{-s-1} \log_\theta\!\lambda\, A (A-\lambda)^{-1} \, x
\, d\lambda = \int_{\Cal C_\theta} \lambda^{-1} \log_\theta\!\lambda\, A
    (A-\lambda)^{-1} \, x \, d\lambda.
\tag A.6
$$
Combining (A.4) and (A.6) (and multiplying with $\frac{i}{2\pi}$) we
obtain the desired result (4.6).

The identity (4.7) stems from \cite{Bu} (we have corrected a sign
here). For this, consider the integration contour
$$
\Gamma_{\theta,\varphi}^N = \{ r e^{i \varphi} \mid N > r > r_0
\} \cup \{ r_0 e^{i\omega} \mid \varphi \ge \omega \ge \theta
\} \cup \{ r e^{i \theta} \mid r_0 < r < N \}.
\tag A.7
$$
Using again $A(A-\lambda)^{-1} = 1 + \lambda(A-\lambda)^{-1}$ we obtain
$$
\int_{\Gamma^N_{\theta,\varphi}} \lambda^{-1} A (A-\lambda)^{-1}\, x
\, d\lambda = \int_{\Gamma^N_{\theta,\varphi}} (A-\lambda)^{-1}\, x \,
d\lambda + \int_{\Gamma^N_{\theta,\varphi}} \lambda^{-1} x \,
d\lambda.
\tag A.8
$$
For the second term we have, using a logarithm with branch
cut disjoint from $\Lambda_{\theta,\varphi}$,
$$
\int_{\Gamma_{\theta,\varphi}^N} \lambda^{-1} d\lambda = \big[\log
\lambda \big]_{N e^{i\varphi}}^{N e^{i\theta}} = i (\theta - \varphi).
\tag A.9
$$
Thus
$$
\frac{i}{2\pi} \int_{\Gamma^N_{\theta,\varphi}} \lambda^{-1} A
(A-\lambda)^{-1}\, x \, d\lambda = \frac{i}{2\pi}
\int_{\Gamma^N_{\theta,\varphi}} (A-\lambda)^{-1}\, x \, d\lambda
+\frac{\varphi-\theta}{2\pi} \, x.  
\tag A.10
$$
For $x\in D(A)$ the limit for $N\to\infty$ is well-defined on
the left-hand side, and the limit of the first term on the right-hand
side then exists as an improper integral, as indicated.
\qed
\enddemo

\demo{Proof of Lemma {\rm 4.3}} The integral along $\Cal C_\theta $ is,
in detail:
$$
\gather
\int_{\Cal C_\theta} \log_\theta \! \lambda \, f(\lambda)\, d\lambda = 
\int_{\infty}^{r_0} (\log r + i \theta) f( r e^{i\theta} ) e^{i\theta}
\,dr + \int_{\theta}^{\theta-2\pi} (\log r_0 + i \omega) f(r_0
e^{i\omega}) i r_0 e^{i\omega} \,d\omega \\ +
\int_{r_0}^{\infty} (\log r + i \theta-2\pi i)) f ( r
e^{i\theta-2\pi i}) e^{i\theta-2\pi i}\, dr.
\tag A.11
\endgather
$$
Since $f(r e^{i \theta - 2\pi i}) e^{i\theta-2\pi i} = f(r e^{i\theta})
e^{i\theta}$, the two terms with $(\log r + i\theta)$ 
cancel each other. Thus
$$
\int_{\Cal C_\theta} \log_\theta \! \lambda \, f(\lambda)\, d\lambda = 
- \int_{\theta-2\pi}^{\theta} (\log r_0 + i \omega)
 f(r_0 e^{i\omega}) i r_0 e^{i\omega}\, d\omega - 2 \pi i \int_{r_0}^{\infty}
f ( r e^{i\theta}) e^{i\theta}\, dr.
\tag A.12
$$
Denote the integrand in the first integral $g(\omega) = (\log r_0 + i
\omega) f(r_0 e^{i\omega}) i r_0 e^{i\omega}$. 

There is of course an identity similar to (A.12) with $\theta$ replaced
by $\varphi$, and then
$$
\aligned
&\int_{\Cal C_{\theta}} \log_\theta \! \lambda \,
f(\lambda) \,d\lambda - \int_{\Cal C_{\varphi}} \log_\varphi \! \lambda \,
f(\lambda)\, d\lambda \\
&= \bigl(-\int_{\theta-2\pi}^\theta  + 
\int_{\varphi-2\pi}^\varphi \bigr)g(\omega) \,d\omega -
2\pi i \big(
\int_{r_0}^\infty f(r e^{i\theta}) e^{i\theta} \,dr 
-\int_{r_0}^\infty f(r e^{i\varphi}) e^{i\varphi} \,dr\big)
\\
&= \bigl(-\int_{\theta-2\pi}^\theta  
+ \int_{\varphi-2\pi}^\varphi \bigr)g(\omega) \,d\omega 
-2\pi i\int_\infty ^{r_0} f(r e^{i\varphi}) e^{i\varphi} \,dr-
2\pi i 
\int_{r_0}^\infty f(r e^{i\theta}) e^{i\theta} \,dr .
\endaligned
\tag A.13
$$
The last two terms are recognized as the contributions to
$-2\pi i\int_{\Gamma _{\theta ,\varphi }}f(\lambda )\,d\lambda $
from the rays 
$e^{i\varphi }[r_0,\infty [\,$ and $e^{i\theta }[r_0,\infty [\,$. 
The first term  is seen to give the contribution from the arc
$\Cal C_{r_0,\theta ,\varphi }=\{r_0e^{i\omega }\mid \varphi 
\ge \omega \ge \theta \}$ as follows: 
$$
\aligned
\bigl(-\int_{\theta-2\pi}^\theta  &+ 
\int_{\varphi-2\pi}^\varphi \bigr)g(\omega) \,d\omega =
\bigl(-\int_{\varphi }^\theta  + 
\int_{\varphi-2\pi}^{\theta -2\pi } \bigr)g(\omega) \,d\omega 
=\int_{\varphi }^\theta [-g(\omega ) + 
g(\omega-2\pi )] \,d\omega \\
&=\int_{\varphi }^\theta [-(\log r_0+i\omega )f(r_0 e^{i\omega}) i r_0
e^{i\omega }+ (\log r_0+i(\omega -2\pi ))f(r_0 e^{i\omega}) i r_0
e^{i\omega }
] \,d\omega \\
&=-2\pi i \int_\varphi ^\theta f(r_0 e^{i\omega}) i r_0 e^{i\omega} \,d\omega
=-2\pi i\int _{\Cal C_{r_0,\theta ,\varphi }}f(\lambda )\,d\lambda .
\quad \square\endaligned
$$
\enddemo

\enddocument